arXiv:1711.09521v1 [math.NA] 27 Nov 2017

# STABILITY AND CONVERGENCE ANALYSIS OF A LINEAR, FULLY DECOUPLED AND UNCONDITIONALLY ENERGY STABLE SCHEME FOR MAGNETO-HYDRODYNAMIC EQUATIONS

GUO-DONG ZHANG*, XIAOMING HE†, AND XIAOFENG YANG‡

**Abstract.** In this paper, we consider numerical approximations for solving the nonlinear magneto-hydrodynamical system, that couples the Navier-Stokes equations and Maxwell equations together. A challenging issue to solve this model numerically is about the time marching problem, i.e., how to develop suitable temporal discretizations for the nonlinear terms in order to preserve the energy stability at the discrete level. We solve this issue in this paper by developing a linear, fully decoupled, first order time stepping scheme, by combining the projection method and some subtle implicit-explicit treatments for nonlinear coupling terms. We further prove that the scheme is unconditional energy stable and derive the optimal error estimates rigorously. Various numerical experiments are implemented to demonstrate the stability and the accuracy in simulating some benchmark simulations, including the Kelvin-Helmholtz shear instability and the magnetic-frozen phenomenon in the lid-driven cavity.

**1. Introduction.** The magneto-hydrodynamical (MHD) system models the behaviors of conducting fluids, such as plasmas, liquid metals, salt water and electrolytes, under external electromagnetic field. It has wide applications in geophysics, astrophysics and confinement for controlled thermonuclear fusion, see [24, 11, 26]. The fundamental concept behind the MHD system is that magnetic fields can induce currents in a moving conductive fluid, which in turn polarizes the fluid and reciprocally changes the magnetic field itself. Thus the governing equations that describe the MHD system are a nonlinear system to couple the Navier-Stokes equations for the hydrodynamics and Maxwell's equations for the electromagnetism. About the extensive theoretical modeling/numerical analysis for the MHD system, we refer to [30, 18, 19, 15, 17, 21, 14, 29, 1, 2, 10, 25, 3, 22] and the references therein.

To solve the MHD system numerically, the main challenging issue is to develop proper temporal discretizations for those coupling terms, including (i) the coupling of between the velocity and pressure in the fluid momentum equation, and (ii) the nonlinear coupling between the magnetic field and the velocity field through convections and Lorentz forces. It is well-known that simple discretizations, like fully explicit or implicit type schemes, can lead to considerable instabilities or suffer from costly time expense. Therefore, people are particularly interested in designing energy stable schemes, in the sense that the discrete energy dissipation laws hold. Meanwhile, while keeping the energy stable feature, it is also desirable to develop schemes that are easy-to-implement. Here the term "easy-to-implement" is referred to "linear" and "decoupled" in comparison with its counter parts: "nonlinear" and "coupled".

It is remarkable that many attempts have been made in this direction recently. In [17], the authors developed two implicit-explicit type methods where the first order method is shown to be unconditionally stable and the second order method is shown to be conditionally stable. However, the model considered in [17] is the reduced version, namely, the magnetic field is assumed to be a fixed function. In [38, 39], the authors developed a decoupled type scheme for the full MHD system, but it is conditionally energy stable with the similar constraint on the time step as well. In [21], the authors developed a totally decoupled scheme where the computations of Navier-Stokes equations are based on the commutator of Laplacian and Leray projection, and all nonlinear and coupling terms are treated explicitly. However, the scheme is still conditionally stable. In [27], the authors had developed a "partially" decoupled scheme where the computations of magnetic field is totally decoupled from the velocity field since all nonlinear terms are treated explicitly, but the velocity is coupled with the pressure in the Navier-Stokes equations. Furthermore, a very severe time step constraint ($\delta t \lesssim h^3$ where $h$ is the grid size of space), which can be very costly in large-scale computations, has to be used to ensure stability. In [7], the authors developed some unconditionally energy stable schemes based on the projection type methods for the Navier-Stokes equations. However, the velocity field and the magnetic field are still coupled together.

Therefore, the first objective of this paper is to develop a time marching scheme that is not only easy-to-implement (*linear and decoupled*), but also unconditionally energy stable. We achieve such a goal by combining several effective approaches, including, (i) an auxiliary intermediate velocity variable to decouple the computation of the magnetic field from the velocity; (ii) the projection method to decouple the pressure from the velocity; and (iii) some subtle implicit-explicit treatments to discretize the nonlinear convection and Lorentz force terms. To the best of the authors' knowledge, this is the first linear, fully decoupled and unconditionally energy stable scheme for solving the MHD system. The second objective is to derive the optimal error estimate for the developed scheme. The final objective is to implement various numerical simulations, including the convergence test, energy stability test and two physical benchmark problems, the Kelvin-Helmholtz shear instability and the magnetic-frozen phenomenon in the lid-driven cavity, to demonstrate the stability and accuracy of the scheme.

*Department of Mathematics, Yantai University, Yantai, 264005, Shandong, P. R. China. Email: gdzhang@ytu.edu.cn. This author's research is partially supported by National Science Foundation of China under grant numbers 11601468 and 11771375.

†Department of Mathematics, Missouri University of Science & Technology, Rolla, MO, 65409, USA. Email: hex@mst.edu. This author's research is partially supported by the U.S. National Science Foundation under grant numbers DMS-1418624 and DMS-1722647.

‡Corresponding author, Department of Mathematics, University of South Carolina, Columbia, SC, 29208, USA. Email: xfyang@math.sc.edu. This author's research is partially supported by the U.S. National Science Foundation under grant numbers DMS-1418898 and DMS-1720212.

The rest of paper is organized as follows. In Section 2, we present the model and derive the associated energy dissipation law. In Section 3, we develop the decoupled scheme and prove its associated energy stability and derive its error estimates. In Section 4, various numerical experiments are presented to demonstrate the stability and effectiveness of the scheme. Finally, some concluding remarks are given in Section 5.

**2. The MHD Model and its energy law.** Here and after, for two vector functions $\boldsymbol{x} = (x_1, x_2), \boldsymbol{y} = (y_1, y_2)$, we denote the $L^2$ inner product as $(\boldsymbol{x}, \boldsymbol{y}) = \int_\Omega \boldsymbol{x}\boldsymbol{y}^T dx$ and $L^2$ norm $\|\boldsymbol{x}\|^2 = (\boldsymbol{x}, \boldsymbol{x})$. We use $H^1(\Omega)$ and $H^2(\Omega)$ to denote the usual sobolev spaces, and use $\|\cdot\|_1$ for the norm in $H^1(\Omega)$ and $\|\cdot\|_2$ for the norm in $H^2(\Omega)$. We also define $H_0^1(\Omega) = \{\phi \in H^1(\Omega) : \phi|_{\partial\Omega} = 0\}$, $L_0^2(\Omega) = \{\phi \in L^2(\Omega) : \int_\Omega \phi dx = 0\}$, and $H(\Omega) = \{\boldsymbol{u} \in L^2(\Omega)^2, \nabla \cdot \boldsymbol{u} = 0, \boldsymbol{u} \cdot \boldsymbol{n}|_{\partial\Omega} = 0\}$.

We consider the following incompressible two dimensional MHD equations:

$$\boldsymbol{u}_t - \nu\Delta\boldsymbol{u} + (\boldsymbol{u} \cdot \nabla)\boldsymbol{u} + \nabla p + s\boldsymbol{B} \times \nabla \times \boldsymbol{B} = 0, \tag{2.1}$$
$$\boldsymbol{B}_t + \eta\nabla \times \nabla \times \boldsymbol{B} - \nabla \times (\boldsymbol{u} \times \boldsymbol{B}) = 0, \tag{2.2}$$
$$\nabla \cdot \boldsymbol{u} = 0, \tag{2.3}$$
$$\nabla \cdot \boldsymbol{B} = 0, \tag{2.4}$$

for $(\boldsymbol{x}, t) \in \Omega \times [0, T)$ with $\Omega \subset \mathbb{R}^2$, where, $\boldsymbol{u}$ denotes the velocity field, $p$ is the pressure, $\boldsymbol{B}$ is the magnetic field. For the physical parameters, $\nu^{-1} = R_e$ (fluid Reynolds number), $\eta^{-1} = R_m$ (magnetic Reynolds number), and $s$ is the coupling coefficient, which are given by

$$R_e = \frac{UL}{\mu_f}, \quad R_m = \mu_m \sigma UL, \quad s = \frac{B^2}{\rho\mu_m U^2},$$

where $U$ is the characteristic velocity, $L$ is the characteristic length, $\mu_f$ is the kinematic viscosity, $\mu_m$ is the magnetic permeability, $\sigma$ is the electric conductivity, $B$ is the characteristic magnetic field, and $\rho$ is the fluid density. The system is equipped with the following boundary conditions

$$\boldsymbol{u}|_{\partial\Omega} = 0, \quad \boldsymbol{B} \times \boldsymbol{n}|_{\partial\Omega} = 0, \tag{2.5}$$

and initial conditions

$$\boldsymbol{u}|_{(t=0)} = \boldsymbol{u}_0(\boldsymbol{x}), \quad \boldsymbol{B}|_{(t=0)} = \boldsymbol{B}_0(\boldsymbol{x}), \tag{2.6}$$

with $\nabla \cdot \boldsymbol{u}_0 = 0$, $\nabla \cdot \boldsymbol{B}_0 = 0$, where $\boldsymbol{n}$ denotes the outward unit normal of $\partial\Omega$.

The model (2.1)-(2.4) follows the energy dissipation law. By taking the $L^2$ innner product of (2.1) with $\boldsymbol{u}$, and of (2.2) with $s\boldsymbol{B}$, using (2.3)-(2.4) and integration by parts, we have

$$(\boldsymbol{u}_t, \boldsymbol{u}) + \nu\|\nabla\boldsymbol{u}\|^2 + s(\boldsymbol{B} \times \nabla \times \boldsymbol{B}, \boldsymbol{u}) = 0, \tag{2.7}$$
$$s(\boldsymbol{B}_t, \boldsymbol{B}) + s\eta\|\nabla \times \boldsymbol{B}\|^2 - s(\nabla \times (\boldsymbol{u} \times \boldsymbol{B}), \boldsymbol{B}) = 0. \tag{2.8}$$

By taking the summation of the two equalities, we obtain

$$\frac{d}{dt}E(\boldsymbol{u}, \boldsymbol{B}) = -\nu\|\nabla\boldsymbol{u}\|^2 - s\eta\|\nabla \times \boldsymbol{B}\|^2, \tag{2.9}$$

where

$$E(\boldsymbol{u}, \boldsymbol{B}) = \frac{1}{2}\|\boldsymbol{u}\|^2 + \frac{s}{2}\|\boldsymbol{B}\|^2 \tag{2.10}$$

represents the total energy of the system (2.1)-(2.4).

## 3. Numerical scheme.

**3.1. Scheme and its Unconditional Energy Stability.** We now construct a semi-discrete time marching numerical scheme for solving the model system (2.1)-(2.4) and prove the corresponding energy stability. It will be clear that the energy stabilities of the semi-discrete schemes are also valid in the fully discrete formulation, for instance by finite element or spectral spatial discretizations.

Let $\delta t > 0$ denote the time step size and set $t_n = n\delta t$ for $0 \le n \le [\frac{T}{\delta t}]$ with the final time $T$. Our numerical scheme reads as follows.

Given the initial conditions $(\boldsymbol{u}^0, \boldsymbol{B}^0, p^0)$, having computed $(\boldsymbol{u}^n, \boldsymbol{B}^n, p^n)$ for $n > 0$, we compute $(\boldsymbol{u}^{n+1}, \boldsymbol{B}^{n+1}, p^{n+1})$ by the following steps.

*Step 1.*

$$\frac{\boldsymbol{B}^{n+1}-\boldsymbol{B}^n}{\delta t}+\eta\nabla\times\nabla\times\boldsymbol{B}^{n+1}-\nabla\times(\boldsymbol{u}^n_\star\times\boldsymbol{B}^n)=0, \tag{3.1}$$

$$\frac{\boldsymbol{u}^n_\star-\boldsymbol{u}^n}{\delta t}+s\boldsymbol{B}^n\times\nabla\times\boldsymbol{B}^{n+1}=0, \tag{3.2}$$

$$\boldsymbol{B}^{n+1}\times\boldsymbol{n}|_{\partial\Omega}=0. \tag{3.3}$$

*Step 2.*

$$\frac{\tilde{\boldsymbol{u}}^{n+1}-\boldsymbol{u}^n_\star}{\delta t}+(\boldsymbol{u}^n\cdot\nabla)\tilde{\boldsymbol{u}}^{n+1}-\nu\Delta\tilde{\boldsymbol{u}}^{n+1}+\nabla p^n=0, \tag{3.4}$$

$$\tilde{\boldsymbol{u}}^{n+1}|_{\partial\Omega}=0. \tag{3.5}$$

*Step 3.*

$$\frac{\boldsymbol{u}^{n+1}-\tilde{\boldsymbol{u}}^{n+1}}{\delta t}+\nabla(p^{n+1}-p^n)=0, \tag{3.6}$$

$$\nabla\cdot\boldsymbol{u}^{n+1}=0, \tag{3.7}$$

$$\boldsymbol{u}^{n+1}\cdot\boldsymbol{n}|_{\partial\Omega}=0. \tag{3.8}$$

A couple remarks are in order.

**Remark 3.1.** *To decouple the computation of the pressure from that of the velocity, we use the first order pressure-correction scheme [12, 31, 35, 13, 41, 42, 34, 43]. To further decouple the computations of $\boldsymbol{B}$ from the velocity field $\boldsymbol{u}$, inspired by [23, 6, 32, 40, 34, 37, 33, 20, 35], we introduce a new, explicit, convective velocity $\boldsymbol{u}^n_\star$, that can be computed directly from (3.2), i.e.,*

$$\boldsymbol{u}^n_\star=\boldsymbol{u}^n+\delta ts\nabla\times\boldsymbol{B}^{n+1}\times\boldsymbol{B}^n. \tag{3.9}$$

*Indeed, if plugging (3.9) into (3.1), one obtain a linear equation for $\boldsymbol{B}^{n+1}$ as*

$$\frac{\boldsymbol{B}^{n+1}}{\delta t}+\eta\nabla\times\nabla\times\boldsymbol{B}^{n+1}+\delta ts\nabla\times\big(\boldsymbol{B}^n\times(\nabla\times\boldsymbol{B}^{n+1}\times\boldsymbol{B}^n)\big)=g^n, \tag{3.10}$$

*where $g^n=-\nabla\times(\boldsymbol{B}^n\times\boldsymbol{u}^n)+\frac{\boldsymbol{B}^n}{\delta t}$. Thus its associated weak form can be written as the following system with the unknown $\boldsymbol{B}\in H^1(\Omega)$,*

$$\frac{1}{\delta t}(\boldsymbol{B},\boldsymbol{C})+\eta(\nabla\times\boldsymbol{B},\nabla\times\boldsymbol{C})+\delta ts(\boldsymbol{B}^n\times\nabla\times\boldsymbol{B},\boldsymbol{B}^n\times\nabla\times\boldsymbol{C})=(g^n,\boldsymbol{C}), \tag{3.11}$$

*for any $\boldsymbol{C}\in H^1(\Omega)$. It is easy to see that the linear operator is symmetric positive definite, and it is quite straight forward to obtain the well-posedness from the Lax-Milgram theorem. Since the proof is rather straight forward, we leave this to the interested readers. Once we obtain $\boldsymbol{B}^{n+1}$ from (3.10), we can update $\tilde{\boldsymbol{u}}^{n+1}$ from*

$$\frac{\tilde{\boldsymbol{u}}^{n+1}-\boldsymbol{u}^n}{\delta t}+s\boldsymbol{B}^n\times\nabla\times\boldsymbol{B}^{n+1}+(\boldsymbol{u}^n\cdot\nabla)\tilde{\boldsymbol{u}}^{n+1}-\nu\Delta\tilde{\boldsymbol{u}}^{n+1}+\nabla p^n=0. \tag{3.12}$$

**Remark 3.2.** *For the pressure equation, indeed, by taking the divergence for (3.6), we get*

$$-\Delta p^{n+1}=-\frac{1}{\delta t}\nabla\cdot\tilde{\boldsymbol{u}}^{n+1}-\Delta p^n, \tag{3.13}$$

*associated with the Neumann boundary conditions $\partial_{\boldsymbol{n}}(p^{n+1}-p^n)|_{\partial\Omega}=0$. Once $p^{n+1}$ is obtained, we update $\boldsymbol{u}^{n+1}$ from $\boldsymbol{u}^{n+1}=\tilde{\boldsymbol{u}}^{n+1}-\delta t\nabla(p^{n+1}-p^n)$.*

The energy stability of the scheme (3.1)-(3.8) is presented as follows.

**Theorem** 3.1. *The scheme (3.1)-(3.8) is unconditionally energy stable in the sense that*

$$\begin{aligned} s\|\boldsymbol{B}^{n+1}\|^2 + \|\boldsymbol{u}^{n+1}\|^2 + \delta t^2\|\nabla p^{n+1}\|^2 + 2\delta t\left(s\eta\|\nabla\times\boldsymbol{B}^{n+1}\|^2 + \nu\|\nabla\tilde{\boldsymbol{u}}^{n+1}\|^2\right) \\ \le s\|\boldsymbol{B}^n\|^2 + \|\boldsymbol{u}^n\|^2 + \delta t^2\|\nabla p^n\|^2. \end{aligned} \tag{3.14}$$

*Moreover,* $\nabla\cdot\boldsymbol{B}^{n+1} = \cdots = \nabla\cdot\boldsymbol{B}^0 = 0$ *for* $1 \le n \le N = [\frac{T}{\delta t}] - 1$.

*Proof.* By taking the $L^2$ inner product of (3.1) with $s\boldsymbol{B}^{n+1}$ and of (3.2) with $\boldsymbol{u}_\star^n$, we obtain

$$\frac{s}{2\delta t}\left(\|\boldsymbol{B}^{n+1}-\boldsymbol{B}^n\|^2 + \|\boldsymbol{B}^{n+1}\|^2 - \|\boldsymbol{B}^n\|^2\right) + s\eta\|\nabla\times\boldsymbol{B}^{n+1}\|^2 + s(\boldsymbol{B}^n\times\boldsymbol{u}_\star^n, \nabla\times\boldsymbol{B}^{n+1}) = 0, \tag{3.15}$$

and

$$\frac{1}{2\delta t}\left(\|\boldsymbol{u}_\star^n - \boldsymbol{u}^n\|^2 + \|\boldsymbol{u}_\star^n\|^2 - \|\boldsymbol{u}^n\|^2\right) - s(\nabla\times\boldsymbol{B}^{n+1}\times\boldsymbol{B}^n, \boldsymbol{u}_\star^n) = 0. \tag{3.16}$$

By taking the $L^2$ inner product of (3.4) with $\tilde{\boldsymbol{u}}^{n+1}$ and using the well known property of

$$(\boldsymbol{u}\cdot\nabla\boldsymbol{v}, \boldsymbol{v}) = 0, \forall\, \boldsymbol{u}\in H, \boldsymbol{v}\in H_0^1(\Omega)^2, \tag{3.17}$$

thus we derive

$$\frac{1}{2\delta t}\left(\|\tilde{\boldsymbol{u}}^{n+1} - \boldsymbol{u}_\star^n\|^2 + \|\tilde{\boldsymbol{u}}^{n+1}\|^2 - \|\boldsymbol{u}_\star^n\|^2\right) + \nu\|\nabla\tilde{\boldsymbol{u}}^{n+1}\|^2 + (\nabla p^n, \tilde{\boldsymbol{u}}^{n+1}) = 0. \tag{3.18}$$

We rewrite (3.6) as

$$\frac{1}{\delta t}\boldsymbol{u}^{n+1} + \nabla p^{n+1} = \frac{1}{\delta t}\tilde{\boldsymbol{u}}^{n+1} + \nabla p^n, \tag{3.19}$$

and take the $L^2$ inner product of the above with itself on both sides, we obtain

$$\frac{1}{2\delta t}\|\boldsymbol{u}^{n+1}\|^2 - \frac{1}{2\delta t}\|\tilde{\boldsymbol{u}}^{n+1}\|^2 + \frac{\delta t}{2}\|\nabla p^{n+1}\|^2 - \frac{\delta t}{2}\|\nabla p^n\|^2 = (\tilde{\boldsymbol{u}}^{n+1}, \nabla p^n). \tag{3.20}$$

Then, by taking the summations (3.15), (3.16), (3.18) and (3.20), we obtain

$$\begin{aligned} &\frac{s}{2\delta t}(\|\boldsymbol{B}^{n+1}\|^2 - \|\boldsymbol{B}^n\|^2 + \|\boldsymbol{B}^{n+1}-\boldsymbol{B}^n\|^2) + s\eta\|\nabla\times\boldsymbol{B}^{n+1}\|^2 \\ &\quad + \frac{1}{2\delta t}\left(\|\boldsymbol{u}^{n+1}\|^2 - \|\boldsymbol{u}^n\|^2 + \|\boldsymbol{u}_\star^n - \boldsymbol{u}^n\|^2\right) + \frac{1}{2\delta t}\|\tilde{\boldsymbol{u}}^{n+1} - \boldsymbol{u}_\star^n\|^2 \\ &\quad + \nu\|\nabla\tilde{\boldsymbol{u}}^{n+1}\|^2 + \frac{\delta t}{2}\|\nabla p^{n+1}\|^2 - \frac{\delta t}{2}\|\nabla p^n\|^2 = 0. \end{aligned} \tag{3.21}$$

After multiplying with $2\delta t$ and dropping some positive terms, we obtain (3.14). Finally, by taking the divergence for (3.1), we get $\nabla\cdot\boldsymbol{B}^{n+1} = \nabla\cdot\boldsymbol{B}^n = \cdots = \nabla\cdot\boldsymbol{B}^0 = 0$. ∎

**3.2. Error analysis.** We now establish the error estimates for scheme (3.1)-(3.8). We denote by $C$ a generic constant that is independent of $\delta t$ but possibly depends on the data and the solution, and use $f \lesssim g$ to say that there is a generic constant $C$ such that $f \le Cg$.

To this end, we rewrite (2.1)-(2.2) as follows.

$$\frac{\boldsymbol{B}(t_{n+1}) - \boldsymbol{B}(t_n)}{\delta t} + \eta\nabla\times\nabla\times\boldsymbol{B}(t_{n+1}) - \nabla\times(\boldsymbol{u}(t_n)\times\boldsymbol{B}(t_n)) = R_b^{n+1}, \tag{3.22}$$

$$\begin{aligned} \frac{\boldsymbol{u}(t_{n+1}) - \boldsymbol{u}(t_n)}{\delta t} - \nu\Delta\boldsymbol{u}(t_{n+1}) + (\boldsymbol{u}(t_n)\cdot\nabla)\boldsymbol{u}(t_{n+1}) + \nabla p(t_n) \\ + s\boldsymbol{B}(t_n)\times\nabla\times\boldsymbol{B}(t_{n+1}) = R_u^{n+1}, \end{aligned} \tag{3.23}$$

$$\frac{\boldsymbol{u}(t_{n+1}) - \boldsymbol{u}(t_{n+1})}{\delta t} + \nabla(p(t_{n+1}) - p(t_n)) = R_p^{n+1}, \tag{3.24}$$

where $R_b^{n+1}, R_u^{n+1}, R_p^{n+1}$ are the truncation errors. We assume the exact solution $(\boldsymbol{u}, \boldsymbol{B}, p)$ of the system (2.1)-(2.4) possesses the following regularity,

$$\text{(3.25)} \qquad \textbf{(A):} \quad \begin{cases} \boldsymbol{u}, \boldsymbol{B} \in L^\infty(0,T;H^2(\Omega)), p \in L^\infty(0,T;H^1(\Omega)), \\ \boldsymbol{u}_t, p_t, \boldsymbol{B}_t \in L^\infty(0,T;H^1(\Omega)), \\ \boldsymbol{u}_{tt}, \boldsymbol{B}_{tt} \in L^\infty(0,T;L^2(\Omega)) \cap L^2(0,T;H^1(\Omega)), p_{tt} \in L^2(0,T;H^1(\Omega)), \\ \boldsymbol{u}_{ttt}, \boldsymbol{B}_{ttt} \in L^2(0,T;L^2(\Omega)). \end{cases}$$

One can easily establish the following estimates for the truncation errors, provided that the exact solutions are sufficiently smooth or satisfy the assumption (A).

**Lemma** 3.1. *Under the Assumption (A), the truncations errors satisfy*

$$\|R_u^n\| + \|R_b^n\| + \|R_p^n\| \lesssim \delta t, 0 \le n \le [\frac{T}{\delta t}]. \tag{3.26}$$

*Proof.* Since the proof is rather standard, due to the page limit, we leave it to the interested readers. □

To derive the error estimates, we denote the error functions as

$$\begin{cases} e_b^n = \boldsymbol{B}(t_n) - \boldsymbol{B}^n, & \tilde{e}_u^n = \boldsymbol{u}(t_n) - \tilde{\boldsymbol{u}}^n, \\ e_u^n = \boldsymbol{u}(t_n) - \boldsymbol{u}^n, & e_p^n = p(t_n) - p^n. \end{cases} \tag{3.27}$$

By subtracting (3.1) from (3.22), (3.4) from (3.23) and applying (3.2), and (3.6) from (3.24), we obtain the following error equations,

$$\frac{e_b^{n+1} - e_b^n}{\delta t} + \eta \nabla \times \nabla \times e_b^{n+1} + \nabla \times (\boldsymbol{B}(t_n) \times \boldsymbol{u}(t_n)) - \nabla \times (\boldsymbol{B}^n \times \boldsymbol{u}_\star^n) = R_b^{n+1}, \tag{3.28}$$

$$\begin{aligned} \frac{\tilde{e}_u^{n+1} - e_u^n}{\delta t} &- \nu \Delta \tilde{e}_u^{n+1} + (\boldsymbol{u}(t_n) \cdot \nabla)\boldsymbol{u}(t_{n+1}) - (\boldsymbol{u}^n \cdot \nabla)\tilde{\boldsymbol{u}}^{n+1} + \nabla e_p^n \\ &+ s\boldsymbol{B}(t_n) \times \nabla \times \boldsymbol{B}(t_{n+1}) - s\boldsymbol{B}^n \times \nabla \times \boldsymbol{B}^{n+1} = R_u^{n+1}, \end{aligned} \tag{3.29}$$

$$\frac{e_u^{n+1} - \tilde{e}_u^{n+1}}{\delta t} + \nabla(e_p^{n+1} - e_p^n) = R_p^{n+1}. \tag{3.30}$$

We first show the $L^\infty$ stability of $\boldsymbol{B}^n$, that plays a key role in the error estimates. Let

$$\kappa = \max_{0 \le t \le T} \|\boldsymbol{B}(t)\|_{L^\infty} + 1, \tag{3.31}$$

the preliminary result is given in the following lemma.

**Lemma** 3.2. *Assuming that the solution to* (2.1)-(2.4) *satisfies Assumption (A), there exists a constant C that is given in the proof, such that, the solution $\boldsymbol{B}^n$ of scheme (3.1)-(3.8) satisfies*

$$\|\boldsymbol{B}^n\|_{L^\infty} \le \kappa, \quad n = 0, 1, \cdots, [\frac{T}{\delta t}], \ \textit{for } \delta t \le C. \tag{3.32}$$

*Proof.* We use mathematical induction to prove this Lemma.

When $n = 0$, we have $\|\boldsymbol{B}^0\|_{L^\infty} \le \kappa$.

Assuming that $\|\boldsymbol{B}^n\|_{L^\infty} \le \kappa$ is valid for $n = 0, 1, \cdots, N$, we show $\|B^{N+1}\|_{L^\infty} \le \kappa$ is also valid through the following three steps.

(i). By taking the $L^2$ inner product of (3.28) with $e_b^{n+1}$ and using integration by parts, we obtain

$$\begin{aligned} \frac{1}{2\delta t}(\|e_b^{n+1} - e_b^n\|^2 + \|e_b^{n+1}\|^2 - \|e_b^n\|^2) + \eta\|\nabla \times e_b^{n+1}\|^2 \\ + (\boldsymbol{B}(t_n) \times \boldsymbol{u}(t_n) - \boldsymbol{B}^n \times \boldsymbol{u}_\star^n, \nabla \times e_b^{n+1}) = (R_b^{n+1}, e_b^{n+1}). \end{aligned} \tag{3.33}$$

By taking the $L^2$ inner product of (3.29) with $\tilde{e}_u^{n+1}$, we derive

$$\begin{aligned} &\frac{1}{2\delta t}(\|\tilde{e}_u^{n+1} - e_u^n\|^2 + \|\tilde{e}_u^{n+1}\|^2 - \|e_u^n\|^2) + \nu\|\nabla \tilde{e}_u^{n+1}\|^2 + (\nabla e_p^n, \tilde{e}_u^{n+1}) \\ &+ \big((\boldsymbol{u}(t_n) \cdot \nabla)\boldsymbol{u}(t_{n+1}) - (\boldsymbol{u}^n \cdot \nabla)\tilde{\boldsymbol{u}}^{n+1}, \tilde{e}_u^{n+1}\big) \\ &+ s\big(\boldsymbol{B}(t_n) \times \nabla \times \boldsymbol{B}(t_{n+1}) - \boldsymbol{B}^n \times \nabla \times \boldsymbol{B}^{n+1}, \tilde{e}_u^{n+1}\big) = (R_u^{n+1}, \tilde{e}_u^{n+1}). \end{aligned} \tag{3.34}$$

We rewrite (3.30) to obtain

$$\frac{1}{\delta t}e_u^{n+1} + \nabla e_p^{n+1} = \frac{1}{\delta t}\tilde{e}_u^{n+1} + \nabla e_p^n + R_p^{n+1}. \tag{3.35}$$

By taking the $L^2$ inner product of (3.35) with itself on both sides, we obtain

$$\begin{aligned}(\tilde{e}_u^{n+1}, \nabla e_p^n) =& \frac{1}{2\delta t}(\|e_u^{n+1}\|^2 - \|\tilde{e}_u^{n+1}\|^2) + \frac{\delta t}{2}(\|\nabla e_p^{n+1}\|^2 - \|\nabla e_p^n\|^2) \\ &- (\tilde{e}_u^{n+1}, R_p^{n+1}) - \delta t(\nabla e_p^n, R_p^{n+1}) - \frac{\delta t}{2}\|R_p^{n+1}\|^2.\end{aligned} \tag{3.36}$$

We combine (3.33)-(3.36) to obtain

$$\begin{aligned}&\frac{1}{2\delta t}(\|e_b^{n+1} - e_b^n\|^2 + \|e_b^{n+1}\|^2 - \|e_b^n\|^2) + \eta\|\nabla \times e_b^{n+1}\|^2 + \frac{\delta t}{2}(\|\nabla e_p^{n+1}\|^2 - \|\nabla e_p^n\|^2) \\ &\quad + \frac{1}{2\delta t}(\|\tilde{e}_u^{n+1} - e_u^n\|^2 + \|e_u^{n+1}\|^2 - \|e_u^n\|^2) + \nu\|\nabla\tilde{e}_u^{n+1}\|^2 \\ &\quad = -(\boldsymbol{B}(t_n) \times \boldsymbol{u}(t_n) - \boldsymbol{B}^n \times \boldsymbol{u}_\star^n, \nabla \times e_b^{n+1}) && (\text{: term A}) \\ &\qquad - \big((\boldsymbol{u}(t_n)\cdot\nabla)\boldsymbol{u}(t_{n+1}) - (\boldsymbol{u}^n\cdot\nabla)\tilde{\boldsymbol{u}}^{n+1}, \tilde{e}_u^{n+1}\big) && (\text{: term B}) \\ &\qquad - s\big(\boldsymbol{B}(t_n) \times \nabla \times \boldsymbol{B}(t_{n+1}) - \boldsymbol{B}^n \times \nabla \times \boldsymbol{B}^{n+1}, \tilde{e}_u^{n+1}\big) && (\text{: term C}) \\ &\qquad + (\tilde{e}_u^{n+1}, R_p^{n+1}) + \delta t(\nabla e_p^n, R_p^{n+1}) && (\text{: term D}) \\ &\qquad + (R_b^{n+1}, e_b^{n+1}) + (R_u^{n+1}, \tilde{e}_u^{n+1}) && (\text{: term E}) \\ &\qquad + \frac{\delta t}{2}\|R_p^{n+1}\|^2.\end{aligned} \tag{3.37}$$

For term A, using the definition of $u_\star^n$ in (3.9) and assumption (A), we derive

$$\begin{aligned}(\text{term A}) &\leq \left|(\boldsymbol{B}(t_n) \times \boldsymbol{u}(t_n) - \boldsymbol{B}^n \times \boldsymbol{u}_\star^n, \nabla \times e_b^{n+1})\right| \\ &= \left|(e_b^n \times \boldsymbol{u}(t_n) + \boldsymbol{B}^n \times (\boldsymbol{u}(t_n) - \boldsymbol{u}_\star^n), \nabla \times e_b^{n+1})\right| \\ &= \left|(e_b^n \times \boldsymbol{u}(t_n), \nabla \times e_b^{n+1}) + (\boldsymbol{B}^n \times e_u^n, \nabla \times e_b^{n+1}) - s\delta t(\boldsymbol{B}^n \times (\nabla \times \boldsymbol{B}^{n+1} \times \boldsymbol{B}^n), \nabla \times e_b^{n+1})\right| \\ &\leq \big(\|e_b^n\|\|\boldsymbol{u}(t_n)\|_{L^\infty} + \|\boldsymbol{B}^n\|_{L^\infty}\|e_u^n\| + s\delta t\|\boldsymbol{B}^n\|_{L^\infty}^2\|\nabla \times \boldsymbol{B}^{n+1}\|\big)\|\nabla \times e_b^{n+1}\| \\ &\lesssim \frac{\eta}{6}\|\nabla \times e_b^{n+1}\|^2 + \|e_b^n\|^2 + \|e_u^n\|^2 + \delta t^2\|\nabla \times \boldsymbol{B}^{n+1}\|^2.\end{aligned} \tag{3.38}$$

For term B, using (3.17), we derive

$$\begin{aligned}(\text{term B}) &\leq \left|((\boldsymbol{u}(t_n)\cdot\nabla)\boldsymbol{u}(t_{n+1}) - (\boldsymbol{u}^n\cdot\nabla)\tilde{\boldsymbol{u}}^{n+1}, \tilde{e}_u^{n+1})\right| \\ &= \left|((e_u^n\cdot\nabla)\boldsymbol{u}(t_{n+1}), \tilde{e}_u^{n+1}) + ((\boldsymbol{u}^n\cdot\nabla)\tilde{e}_u^{n+1}, \tilde{e}_u^{n+1})\right| \\ &= \left|((e_u^n\cdot\nabla)\boldsymbol{u}(t_{n+1}), \tilde{e}_u^{n+1})\right| = \left|((e_u^n\cdot\nabla)\tilde{e}_u^{n+1}, \boldsymbol{u}(t_{n+1}))\right| \\ &\leq \|e_u^n\|\,\|\nabla\tilde{e}_u^{n+1}\|\|\boldsymbol{u}(t_{n+1})\|_{L^\infty} \lesssim \frac{\nu}{8}\|\nabla\tilde{e}_u^{n+1}\|^2 + \|e_u^n\|^2.\end{aligned} \tag{3.39}$$

For term C, using $\tilde{e}_u^{n+1} = e_u^{n+1} + \delta t(\nabla e_p^{n+1} - \nabla e_p^n) - \delta t R_p^{n+1}$, we derive

$$\begin{aligned}(\text{term C}) &\leq s\left|(\nabla \times \boldsymbol{B}(t_{n+1}) \times \boldsymbol{B}(t_n) - \nabla \times \boldsymbol{B}^{n+1} \times \boldsymbol{B}^n, \tilde{e}_u^{n+1})\right| \\ &= s\left|(\nabla \times \boldsymbol{B}(t_{n+1}) \times e_b^n, \tilde{e}_u^{n+1}) + (\nabla \times e_b^{n+1} \times \boldsymbol{B}^n, \tilde{e}_u^{n+1})\right| \\ &\leq s\|\nabla \times \boldsymbol{B}(t_{n+1})\|_{L^4}\|e_b^n\|\,\|\tilde{e}_u^{n+1}\|_{L^4} + s\|\nabla \times e_b^{n+1}\|\,\|\boldsymbol{B}^n\|_{L^\infty}\|\tilde{e}_u^{n+1}\| \\ &\lesssim \|e_b^n\|\|\nabla\tilde{e}_u^{n+1}\| + \|\nabla \times e_b^{n+1}\|(\|e_u^{n+1}\| + \delta t\|\nabla e_p^{n+1} - \nabla e_p^n\| + \delta t\|R_p^{n+1}\|) \\ &\lesssim \frac{\nu}{8}\|\nabla\tilde{e}_u^{n+1}\|^2 + \frac{\eta}{6}\|\nabla \times e_b^{n+1}\|^2 + \|e_b^n\|^2 + \|e_u^{n+1}\|^2 + \delta t^2\|\nabla e_p^{n+1} - \nabla e_p^n\|^2 + \delta t^2\|R_p^{n+1}\|^2.\end{aligned} \tag{3.40}$$

For term D and term E, we derive

$$
\begin{aligned}
\text{(term D)} &\leq |(\tilde{e}_u^{n+1}, R_p^{n+1})| + \delta t|(\nabla e_p^n, R_p^{n+1})| \\
&\leq \|\tilde{e}_u^{n+1}\|\|R_p^{n+1}\| + \delta t\|\nabla e_p^n\|\|R_p^{n+1}\| \\
&\lesssim \frac{\nu}{8}\|\nabla \tilde{e}_u^{n+1}\|^2 + \|R_p^{n+1}\|^2 + \delta t^2\|\nabla e_p^n\|^2 + \|R_p^{n+1}\|^2,
\end{aligned} \tag{3.41}
$$

and

$$
\begin{aligned}
\text{(term E)} &\leq |(R_b^{n+1}, e_b^{n+1})| + |(R_u^{n+1}, \tilde{e}_u^{n+1})| \\
&\leq \|R_b^{n+1}\|\ \|e_b^{n+1}\| + \|R_u^{n+1}\|\ \|\tilde{e}_u^{n+1}\| \\
&\lesssim \|R_b^{n+1}\|\ \|\nabla \times e_b^{n+1}\| + \|R_u^{n+1}\|\ \|\nabla \tilde{e}_u^{n+1}\| \\
&\lesssim \frac{\eta}{6}\|\nabla \times e_b^{n+1}\|^2 + \frac{\nu}{8}\|\nabla \tilde{e}_u^{n+1}\|^2 + \|R_b^{n+1}\|^2 + \|R_u^{n+1}\|^2.
\end{aligned} \tag{3.42}
$$

By combining the above estimates with (3.37), we obtain

$$
\begin{aligned}
&\|e_b^{n+1}\|^2 - \|e_b^n\|^2 + \delta t\eta\|\nabla \times e_b^{n+1}\|^2 + \|e_u^{n+1}\|^2 - \|e_u^n\|^2 + \delta t\nu\|\nabla \tilde{e}_u^{n+1}\|^2 + \delta t^2(\|\nabla e_p^{n+1}\|^2 - \|\nabla e_p^n\|^2) \\
&\quad \lesssim \delta t\|e_b^n\|^2 + \delta t\|e_u^n\|^2 + \delta t\|e_u^{n+1}\|^2 + \delta t^3\|\nabla e_p^{n+1} - \nabla e_p^n\|^2 + \delta t^3\|\nabla e_p^n\|^2 \\
&\quad + \delta t^3\|\nabla \times \boldsymbol{B}^{n+1}\|^2 + \delta t^3\|R_p^{n+1}\|^2 + \delta t\|R_p^{n+1}\|^2 \\
&\quad + \delta t\|R_b^{n+1}\|^2 + \delta t\|R_u^{n+1}\|^2 + \delta t^2\|R_p^{n+1}\|^2.
\end{aligned} \tag{3.43}
$$

Summing up the above inequality from $n = 0$ to $m \leq N$ and using Lemma 3.1, we obtain

$$
\begin{aligned}
&\|e_b^{m+1}\|^2 + \|e_u^{m+1}\|^2 + \delta t^2\|\nabla e_p^{m+1}\|^2 + \delta t\sum_{n=0}^{m}\left(\eta\|\nabla \times e_b^{n+1}\|^2 + \nu\|\nabla \tilde{e}_u^{n+1}\|^2\right) \\
&\qquad \lesssim \delta t\sum_{n=0}^{m}\left(\|e_b^{n+1}\|^2 + \|e_u^{n+1}\|^2 + \delta t^2\|\nabla e_p^{n+1}\|^2\right) + \delta t^2,
\end{aligned} \tag{3.44}
$$

where we use $\delta t^3 \sum_{n=0}^{m} \|\nabla \times \boldsymbol{B}^{n+1}\|^2 \lesssim \delta t^2$ that is obtained from Theorem 3.1.

Therefore, using the discrete Gronwall's inequality, there exists positive constants $C_0$ and $C_1$, such that

$$
\|e_b^{m+1}\|^2 + \|e_u^{m+1}\|^2 + \delta t^2\|\nabla e_p^{m+1}\|^2 + \delta t\sum_{n=0}^{m}\left(\eta\|\nabla \times e_b^{n+1}\|^2 + \nu\|\nabla \tilde{e}_u^{n+1}\|^2\right) \leq C_1\delta t^2, \tag{3.45}
$$

for $\delta t \leq C_0$ and $m \leq N$.

(ii). For $n \leq N$, for (3.1), by the $H^2$ regularity of the elliptic problem, there exists a constant $C_\Omega$ that only depends on $\Omega$, such that

$$
\begin{aligned}
\|\boldsymbol{B}^{n+1}\|_2 &\leq C_\Omega\left\|\frac{\boldsymbol{B}^{n+1} - \boldsymbol{B}^n}{\delta t}\right\| + C_\Omega\|\nabla \times (\boldsymbol{B}^n \times \boldsymbol{u}_\star^n)\| \\
&\lesssim \left\|\frac{\boldsymbol{B}^{n+1} - \boldsymbol{B}^n}{\delta t}\right\| + \Big(\|\boldsymbol{u}_\star^n\nabla \boldsymbol{B}^n\| + \|\boldsymbol{B}^n\nabla \boldsymbol{u}_\star^n\| + \|\boldsymbol{B}^n\nabla \cdot \boldsymbol{u}_\star^n\|\Big),
\end{aligned} \tag{3.46}
$$

where we use $\nabla \cdot \boldsymbol{B}^n = 0$ and the following identity

$$
\nabla \times (\boldsymbol{a} \times \boldsymbol{b}) = \boldsymbol{b} \cdot \nabla \boldsymbol{a} - \boldsymbol{a} \cdot \nabla \boldsymbol{b} + \boldsymbol{a}\nabla \cdot \boldsymbol{b} - \boldsymbol{b}\nabla \cdot \boldsymbol{a}. \tag{3.47}
$$

For the first term of (3.46) on the right hand side, from (3.45) and assumption (A), there exists a constant $C_2$ such that

$$
\begin{aligned}
\left\|\frac{\boldsymbol{B}^{n+1} - \boldsymbol{B}^n}{\delta t}\right\| &= \left\|\frac{-e_b^{n+1} + e_b^n}{\delta t} + \frac{\boldsymbol{B}(t_{n+1}) - \boldsymbol{B}(t_n)}{\delta t}\right\| \\
&\leq \left\|\frac{e_b^{n+1} - e_b^n}{\delta t}\right\| + \left\|\frac{\boldsymbol{B}(t_{n+1}) - \boldsymbol{B}(t_n)}{\delta t}\right\| \leq C_2.
\end{aligned} \tag{3.48}
$$

For other terms of (3.46) on the right hand side, by combining (3.2) and (3.6) at $t = t^n$, we obtain

$$\boldsymbol{u}_\star^n = \tilde{\boldsymbol{u}}^n - \delta t \nabla(p^n - p^{n-1}) + \delta t s \nabla \times \boldsymbol{B}^{n+1} \times \boldsymbol{B}^n, \tag{3.49}$$

thus we derive

$$\begin{aligned}
\|\boldsymbol{u}_\star^n \nabla \boldsymbol{B}^n\| &\leq \|\tilde{\boldsymbol{u}}^n \nabla \boldsymbol{B}^n\| + \delta t \|\nabla(p^n - p^{n-1}) \nabla \boldsymbol{B}^n\| + \delta t s \|(\nabla \times \boldsymbol{B}^{n+1} \times \boldsymbol{B}^n) \nabla \boldsymbol{B}^n\| \\
&\leq \|\tilde{e}_u^n \nabla \boldsymbol{B}^n\| + \|\boldsymbol{u}(t_n) \nabla \boldsymbol{B}^n\| + \delta t \|\nabla(p^n - p^{n-1}) \nabla \boldsymbol{B}^n\| \\
&\quad + \delta t s \|(\nabla \times \boldsymbol{B}^{n+1} \times \boldsymbol{B}^n) \nabla \boldsymbol{B}^n\| \\
&\lesssim \|\tilde{e}_u^n\|_{L^4} \|\nabla \boldsymbol{B}^n\|_{L^4} + \|\boldsymbol{u}(t_n)\|_{L^\infty} \|\nabla \boldsymbol{B}^n\| + \delta t \|\nabla(p^n - p^{n-1})\|_{L^4} \|\nabla \boldsymbol{B}^n\|_{L^4} \\
&\quad + \delta t \|\nabla \times \boldsymbol{B}^{n+1}\|_{L^4} \|\boldsymbol{B}^n\|_{L^\infty} \|\nabla \boldsymbol{B}^n\|_{L^4} \\
&\lesssim \|\nabla \tilde{e}_u^n\| \|\nabla \boldsymbol{B}^n\|^{\frac{1}{2}} \|\boldsymbol{B}^n\|_2^{\frac{1}{2}} + \|\nabla \boldsymbol{B}^n\| \\
&\quad + \delta t \|\nabla(p^n - p^{n-1})\|^{\frac{1}{2}} \|p^n - p^{n-1}\|_2^{\frac{1}{2}} \|\nabla \boldsymbol{B}^n\|^{\frac{1}{2}} \|\boldsymbol{B}^n\|_2^{\frac{1}{2}} \\
&\quad + \delta t \|\nabla \times \boldsymbol{B}^{n+1}\|^{\frac{1}{2}} \|\boldsymbol{B}^{n+1}\|_2^{\frac{1}{2}} \|\boldsymbol{B}^n\|_{L^\infty} \|\nabla \boldsymbol{B}^n\|^{\frac{1}{2}} \|\boldsymbol{B}^n\|_2^{\frac{1}{2}} \\
&\lesssim \|\nabla \boldsymbol{B}^n\| + \delta t \|\boldsymbol{B}^n\|_2 + \|\nabla(p^n - p^{n-1})\| \|\nabla \boldsymbol{B}^n\| + \delta t^2 \|p^n - p^{n-1}\|_2 \|\boldsymbol{B}^n\|_2 \\
&\quad + \frac{1}{6} \|\boldsymbol{B}^{n+1}\|_2 + \delta t^2 \|\nabla \times \boldsymbol{B}^{n+1}\| \|\boldsymbol{B}^n\|_{L^\infty}^2 \|\nabla \boldsymbol{B}^n\| \|\boldsymbol{B}^n\|_2 \\
&\lesssim C_3 + \frac{1}{6} \|\boldsymbol{B}^{n+1}\|_2 + \delta t \|\boldsymbol{B}^n\|_2 + \delta t^2 \|\boldsymbol{B}^n\|_2 + \delta t^2 \|p^n - p^{n-1}\|_2 \|\boldsymbol{B}^n\|_2,
\end{aligned} \tag{3.50}$$

where we actually use $\|\nabla(p^n - p^{n-1})\|$, $\|\nabla \boldsymbol{B}^n\|$ and $\|\nabla \times \boldsymbol{B}^{n+1}\|$ are all bounded (this can simply proved from $\|\nabla p^n\| \leq \|\nabla e_p^n\| + \|\nabla p(t^n)\|$, $\|\nabla \boldsymbol{B}^n\| \leq \|\nabla \times e_b^n\| + \|\nabla \boldsymbol{B}(t^n)\|$ by using (3.45) and assumption (A)), and use $\|\nabla \tilde{e}_u^n\| \lesssim \delta t^{\frac{1}{2}}$ that is obtained from (3.45).

Similarly, we derive

$$\begin{aligned}
\|\boldsymbol{B}^n \nabla \boldsymbol{u}_\star^n\| &\leq \|\boldsymbol{B}^n \nabla \tilde{\boldsymbol{u}}^n\| + \delta t \|\boldsymbol{B}^n \nabla \nabla(p^n - p^{n-1})\| + \delta t s \|\boldsymbol{B}^n \nabla(\nabla \times \boldsymbol{B}^{n+1} \times \boldsymbol{B}^n)\| \\
&\lesssim \|\boldsymbol{B}^n\|_{L^\infty} \|\nabla \tilde{\boldsymbol{u}}^n\| + \delta t \|\boldsymbol{B}^n\|_{L^\infty} \|p^n - p^{n-1}\|_2 \\
&\quad + \delta t (\|\boldsymbol{B}^n\|_{L^\infty}^2 \|\boldsymbol{B}^{n+1}\|_2 + \|\boldsymbol{B}^n\|_{L^\infty} \|\nabla \boldsymbol{B}^{n+1}\|_{L^4} \|\nabla \boldsymbol{B}^n\|_{L^4}) \\
&\lesssim \|\boldsymbol{B}^n\|_{L^\infty} \|\nabla \tilde{\boldsymbol{u}}^n\| + \delta t \|\boldsymbol{B}^n\|_{L^\infty} \|p^n - p^{n-1}\|_2 \\
&\quad + \delta t (\|\boldsymbol{B}^n\|_{L^\infty}^2 \|\boldsymbol{B}^{n+1}\|_2 + \|\boldsymbol{B}^n\|_{L^\infty} \|\nabla \boldsymbol{B}^{n+1}\|^{\frac{1}{2}} \|\boldsymbol{B}^{n+1}\|_2^{\frac{1}{2}} \|\nabla \boldsymbol{B}^n\|^{\frac{1}{2}} \|\boldsymbol{B}^n\|_2^{\frac{1}{2}}) \\
&\lesssim \|\boldsymbol{B}^n\|_{L^\infty} \|\nabla \tilde{\boldsymbol{u}}^n\| + \delta t \|\boldsymbol{B}^n\|_{L^\infty} \|p^n - p^{n-1}\|_2 + \delta t \|\boldsymbol{B}^n\|_{L^\infty}^2 \|\boldsymbol{B}^{n+1}\|_2 \\
&\quad + \frac{1}{6} \|\boldsymbol{B}^{n+1}\|_2 + \delta t^2 \|\boldsymbol{B}^n\|_{L^\infty}^2 \|\nabla \boldsymbol{B}^{n+1}\| \|\nabla \boldsymbol{B}^n\| \, \|\boldsymbol{B}^n\|_2 \\
&\lesssim C_4 + \frac{1}{6} \|\boldsymbol{B}^{n+1}\|_2 + \delta t \|p^n - p^{n-1}\|_2 + \delta t^2 \|\boldsymbol{B}^n\|_2 + \delta t \|\boldsymbol{B}^{n+1}\|_2,
\end{aligned} \tag{3.51}$$

where we use $\|\nabla \tilde{\boldsymbol{u}}^n\|$ and $\|\nabla \boldsymbol{B}^n\|$ is bounded by constants.

Likewise, for the last term in (3.46), we have

$$\begin{aligned}
\|\boldsymbol{B}^n \nabla \cdot \boldsymbol{u}_\star^n\| &\leq \|\boldsymbol{B}^n \nabla \cdot \tilde{\boldsymbol{u}}^n\| + \delta t \|\boldsymbol{B}^n \Delta(p^n - p^{n-1})\| + \delta t s \|\boldsymbol{B}^n \nabla \cdot (\nabla \times \boldsymbol{B}^{n+1} \times \boldsymbol{B}^n)\| \\
&\lesssim \|\boldsymbol{B}^n\|_{L^\infty} \|\nabla \cdot \tilde{\boldsymbol{u}}^n\| + \delta t \|\boldsymbol{B}^n\|_{L^\infty} \|p^n - p^{n-1}\|_2 + \delta t \|\boldsymbol{B}^n\|_{L^\infty}^2 \|\boldsymbol{B}^{n+1}\|_2 + \delta t \|\boldsymbol{B}^n\|_{L^\infty} \|\nabla \boldsymbol{B}^{n+1}\|_{L^4} \|\nabla \boldsymbol{B}^n\|_{L^4} \\
&\lesssim C_5 + \frac{1}{6} \|\boldsymbol{B}^{n+1}\|_2 + \delta t \|p^n - p^{n-1}\|_2 + \delta t^2 \|\boldsymbol{B}^n\|_2 + \delta t \|\boldsymbol{B}^{n+1}\|_2.
\end{aligned} \tag{3.52}$$

By taking the divergence for (3.6), we obtain

$$-\delta t \Delta(p^n - p^{n-1}) = -\nabla \cdot \tilde{\boldsymbol{u}}^n.$$

From (3.45), we find

$$\delta t \|p^n - p^{n-1}\|_2 \lesssim \|\nabla \cdot \tilde{\boldsymbol{u}}^n\| = \|\nabla \cdot \tilde{e}_u^n\| \lesssim C_6. \tag{3.53}$$

By combining (3.46)–(3.52) and using (3.53), if $\delta t \leq \widehat{C}_0$, there exist two positive constants $C_7$ and $C_8$ such that

$$\|\boldsymbol{B}^{n+1}\|_2 \leq C_7 + C_8\delta t\|\boldsymbol{B}^n\|_2. \tag{3.54}$$

Therefore, if $C_8\delta t \leq \frac{1}{2}$, i.e., $\delta t \leq \frac{1}{2C_8}$, we have

$$\begin{aligned}\|\boldsymbol{B}^{N+1}\|_2 &\leq \max\{C_7, \|\boldsymbol{B}_0\|_2\}\left(1 + C_8\delta t + (C_8\delta t)^2 + \cdots + (C_8\delta t)^{N+1}\right)\\ &= \max\{C_7, \|\boldsymbol{B}_0\|_2\}\frac{1-(C_8\delta t)^{N+2}}{1-C_8\delta t} \leq 2\max\{C_7, \|\boldsymbol{B}_0\|_2\}.\end{aligned} \tag{3.55}$$

(iii). From (3.55) and the assumption (A), there exists a positive constant $C_9$ such that

$$\|e_b^{N+1}\|_2 \leq \|\boldsymbol{B}^{N+1}\|_2 + \|\boldsymbol{B}(t_{N+1})\|_2 \leq C_9. \tag{3.56}$$

Finally, from (3.45), we have

$$\begin{aligned}\|\boldsymbol{B}^{N+1}\|_{L^\infty} &\leq \|e_b^{N+1}\|_{L^\infty} + \|\boldsymbol{B}(t_{N+1})\|_{L^\infty}\\ &\leq C_{10}\|e_b^{N+1}\|_2^{\frac{1}{2}}\|e_b^{N+1}\|^{\frac{1}{2}} + \|\boldsymbol{B}(t_{N+1})\|_{L^\infty}\\ &\leq C_{10}\sqrt{C_9}\sqrt[4]{C_1}\sqrt{\delta t} + \|\boldsymbol{B}(t_{N+1})\|_{L^\infty}.\end{aligned} \tag{3.57}$$

Thus, if $C_{10}\sqrt{C_9}\sqrt[4]{C_1}\sqrt{\delta t} \leq 1$, i.e., $\delta t \leq \frac{1}{C_{10}^2 C_9\sqrt{C_1}}$, we have

$$\|\boldsymbol{B}^{N+1}\|_{L^\infty} \leq 1 + \|\boldsymbol{B}(t_{N+1})\|_{L^\infty} \leq \kappa. \tag{3.58}$$

Then we obtain the conclusion (3.32) by induction for $\delta t \leq C = \min(C_0, \widehat{C}_0, \frac{1}{2C_8}, \frac{1}{C_{10}^2 C_9\sqrt{C_1}})$. ◻

**Theorem** 3.2. *Suppose the solution to (2.1)-(2.4) satisfies Assumption (A). Then, the scheme (3.1)-(3.6) is unconditionally convergent and has the following error estimates: for* $0 \leq m \leq [\frac{T}{\delta t}] - 1$,

$$\|e_b^{m+1}\|^2 + \|e_u^{m+1}\|^2 + \delta t^2\|\nabla e_p^{m+1}\|^2 + \delta t\sum_{n=0}^{m}\left(\eta\|\nabla\times e_b^{n+1}\|^2 + \nu\|\nabla\tilde{e}_u^{n+1}\|^2\right) \lesssim \delta t^2.$$

*Proof.* Since $\|\boldsymbol{B}^n\|_{L^\infty} \leq \kappa$ for any $0 \leq n \leq [\frac{T}{\delta t}]$, by following the first step of the proof of Lemma 3.2, we obtain that (3.45) is valid for any $0 \leq m \leq [\frac{T}{\delta t}] - 1$ provided $\delta t \leq C$. Besides, if $\delta t \geq C$, using (3.14) and assumption (A), we deduce that there exists a constant $C_{11}$ such that

$$\begin{aligned}&\|e_b^{m+1}\|^2 + \|e_u^{m+1}\|^2 + \delta t^2\|\nabla e_p^{m+1}\|^2\\ &\quad + \delta t\sum_{n=0}^{m}\left(\eta\|\nabla\times e_b^{n+1}\|^2 + \nu\|\nabla\tilde{e}_u^{n+1}\|^2\right) \leq C_{11} = \frac{C_{11}}{C^2}C^2 \leq \frac{C_{11}}{C^2}(\delta t)^2 \lesssim \delta t^2.\end{aligned}$$

Therefore, the proof is finished by combining the two cases. ◻

**Remark** 3.1. *Once the bound* $\|\boldsymbol{B}^n\|_{L^\infty} \leq \kappa$ *and the error estimates in Theorem 3.2 are obtained, by following the step 2 of the proof for Lemma 3.2, we can establish the* $H^2$ *stability of* $\boldsymbol{B}^n$*, namely, there exist a constant* $\widehat{\kappa}$*, such that*

$$\max_{0\leq n\leq[\frac{T}{\delta t}]}\|\boldsymbol{B}^n\|_2 \leq \widehat{\kappa}.$$

**Remark** 3.2. *Define* $P_H$ *as the* $L^2$ *orthogonal projector from* $L^2(\Omega)^2$ *to* $H$*, i.e.,*

$$(\boldsymbol{u} - P_H\boldsymbol{u}, \boldsymbol{v}) = 0, \quad \forall\boldsymbol{u} \in L^2(\Omega)^d, \quad \boldsymbol{v} \in H.$$

*By the* $H^1$ *stability of* $P_H$ *[36], we have*

$$\|e_u^n\|_1 = \|P_H\tilde{e}_u^n\|_1 \leq \|\tilde{e}_{\boldsymbol{u}}^n\|_1, \tag{3.59}$$

*which implies*

$$\delta t \sum_{n=0}^{N} \|\nabla e_u^n\|^2 \lesssim \delta t^2. \tag{3.60}$$

Note that the order for the the pressure is not optimal, therefore we need improve it by the following process. We denote $d_t w^n = \frac{w^n - w^{n-1}}{\delta t}$ for any variable $w$. By applying $d_t$ to (3.28)-(3.30), we obtain

$$\begin{aligned}
&\frac{d_t\tilde{e}_u^{n+1} - d_t e_u^n}{\delta t} - \nu\Delta d_t\tilde{e}_u^{n+1} + \boldsymbol{u}(t_n)\nabla d_t\boldsymbol{u}(t_{n+1}) + d_t\boldsymbol{u}(t_n)\nabla\boldsymbol{u}(t_n) - \boldsymbol{u}^n\nabla d_t\tilde{\boldsymbol{u}}^{n+1} - d_t\boldsymbol{u}^n\nabla\tilde{\boldsymbol{u}}^n \\
&\quad +\nabla d_t e_p^n + s\boldsymbol{B}(t_n)\times\nabla\times d_t\boldsymbol{B}(t_{n+1}) + s d_t\boldsymbol{B}(t_n)\times\nabla\times\boldsymbol{B}(t_n) \\
&\quad -s\boldsymbol{B}^n\times\nabla\times d_t\boldsymbol{B}^{n+1} - s d_t\boldsymbol{B}^n\times\nabla\times\boldsymbol{B}^n = d_t R_u^{n+1},
\end{aligned} \tag{3.61}$$

$$\frac{d_t e_u^{n+1} - d_t\tilde{e}_u^{n+1}}{\delta t} + \nabla(d_t e_p^{n+1} - d_t e_p^n) = d_t R_p^{n+1}, \tag{3.62}$$

$$\begin{aligned}
&\frac{d_t e_b^{n+1} - d_t e_b^n}{\delta t} + \eta\nabla\times\nabla\times d_t e_b^{n+1} + \nabla\times(d_t\boldsymbol{B}(t_n)\times\boldsymbol{u}(t_n)) + \nabla\times(\boldsymbol{B}(t_{n-1})\times d_t\boldsymbol{u}(t_n)) \\
&\quad -\nabla\times(d_t\boldsymbol{B}^n\times\boldsymbol{u}_\star^n) - \nabla\times(\boldsymbol{B}^{n-1}\times d_t\boldsymbol{u}_\star^n) = d_t R_b^{n+1}.
\end{aligned} \tag{3.63}$$

**Lemma** 3.3. *Under the Assumption (A), the truncations errors satisfy*

$$\delta t \sum_{n=0}^{[\frac{T}{\delta t}]} (\|d_t R_u^n\|^2 + \|d_t R_b^n\|^2 + \|d_t R_p^n\|^2) \lesssim \delta t^2. \tag{3.64}$$

*Proof.* Since the proof is rather standard, due to the page limit, we leave it to the interested readers. ∎

**Lemma** 3.4. *Under the assumption (A), there holds*

$$\|d_t e_u^1\|^2 + \|d_t e_b^1\|^2 + \delta t^2\|\nabla d_t e_p^1\|^2 \lesssim \delta t^2. \tag{3.65}$$

*Proof.* (i). by taking $n = 0$ in (3.28) and from $e_b^0 = 0$, we obtain

$$\frac{e_b^1}{\delta t} + \eta\nabla\times\nabla\times e_b^1 + \nabla\times(\boldsymbol{B}_0\times\boldsymbol{u}_0) - \nabla\times(\boldsymbol{B}_0\times\boldsymbol{u}_\star^0) = R_b^1.$$

By taking the $L^2$ inner product of above equation with $\frac{1}{\delta t}e_b^1$, we have

$$\|\frac{e_b^1}{\delta t}\|^2 + \frac{\eta}{\delta t}\|\nabla\times e_b^1\|^2 + (\nabla\times(\boldsymbol{B}_0\times\boldsymbol{u}_0) - \nabla\times(\boldsymbol{B}_0\times\boldsymbol{u}_\star^0), \frac{e_b^1}{\delta t}) = (R_b^1, \frac{e_b^1}{\delta t}). \tag{3.66}$$

From (3.9), we obtain

$$\begin{aligned}
&\left|(\nabla\times(\boldsymbol{B}_0\times\boldsymbol{u}_0) - \nabla\times(\boldsymbol{B}_0\times\boldsymbol{u}_\star^0), \frac{e_b^1}{\delta t})\right| \\
&\quad = \delta t s\left|\left(\nabla\times(\boldsymbol{B}_0\times(\boldsymbol{B}_0\times\nabla\times\boldsymbol{B}^1)), \frac{e_b^1}{\delta t}\right)\right| \\
&\quad = \delta t s\left|\left((\boldsymbol{B}_0\times\nabla\times\boldsymbol{B}^1)\nabla\boldsymbol{B}_0 - \boldsymbol{B}_0\nabla(\boldsymbol{B}_0\times\nabla\times\boldsymbol{B}^1) + \boldsymbol{B}_0\nabla\cdot(\boldsymbol{B}_0\times\nabla\times\boldsymbol{B}^1), \frac{e_b^1}{\delta t}\right)\right| \\
&\quad \lesssim \delta t\|\boldsymbol{B}_0\|_{L^\infty}\|\boldsymbol{B}_0\|_2\|\boldsymbol{B}^1\|_2\|\frac{e_b^1}{\delta t}\| + \delta t\|\boldsymbol{B}_0\|_{L^\infty}\|\boldsymbol{B}_0\|_2\|\boldsymbol{B}^1\|_2\|\frac{e_b^1}{\delta t}\| + \delta t\|\boldsymbol{B}^0\|_{L^\infty}^2\|\boldsymbol{B}^1\|_2\|\frac{e_b^1}{\delta t}\| \\
&\quad \lesssim \frac{1}{4}\|\frac{e_b^1}{\delta t}\|^2 + \delta t^2,
\end{aligned} \tag{3.67}$$

and

$$|(R_b^1, \frac{e_b^1}{\delta t})| \le \|R_b^1\|\|\frac{e_b^1}{\delta t}\| \lesssim \frac{1}{4}\|\frac{e_b^1}{\delta t}\|^2 + \|R_b^1\|^2 \lesssim \frac{1}{4}\|\frac{e_b^1}{\delta t}\|^2 + \delta t^2,$$

where we use $\|\boldsymbol{B}^1\|_2$ is bounded by a constant from Remark 3.1.

Therefore, since $e_b^0 = 0$, (3.66) implies

$$\|\frac{e_b^1 - e_b^0}{\delta t}\|^2 + \frac{\eta}{\delta t}\|\nabla \times e_b^1\|^2 \lesssim \delta t^2. \tag{3.68}$$

(ii). By taking $n = 0$ in (3.29) and using $e_u^0 = e_p^0 = 0$, we obtain

$$\frac{\tilde{e}_u^1}{\delta t} - \nu\Delta\tilde{e}_u^1 + (\boldsymbol{u}_0 \cdot \nabla)\tilde{e}_u^1 + s\boldsymbol{B}_0 \times \nabla \times e_b^1 = R_u^1.$$

By taking $L^2$ inner product of the above equation with $\frac{1}{\delta t}\tilde{e}_u^1$, we have

$$\|\frac{\tilde{e}_u^1}{\delta t}\|^2 + \frac{1}{\delta t}\nu\|\nabla\tilde{e}_u^1\|^2 + s(\boldsymbol{B}_0 \times \nabla \times e_b^1, \frac{\tilde{e}_u^1}{\delta t}) = (R_u^1, \frac{\tilde{e}_u^1}{\delta t}). \tag{3.69}$$

By taking $n = 0$ in (3.30), it gives

$$\frac{e_u^1}{\delta t} + \nabla e_p^1 = \frac{\tilde{e}_u^1}{\delta t} + R_p^1.$$

By taking the $L^2$ inner products of the above equation with itself on both sides, we obtain

$$\frac{1}{2}(\|\frac{e_u^1}{\delta t}\|^2 - \|\frac{\tilde{e}_u^1}{\delta t}\|^2) + \frac{1}{2}\|\nabla e_p^1\|^2 - (\frac{\tilde{e}_u^1}{\delta t}, R_p^1) - \frac{1}{2}\|R_p^1\|^2 = 0. \tag{3.70}$$

We combine (3.69) and (3.70) to obtain

$$\begin{aligned}\frac{1}{2}\|\frac{\tilde{e}_u^1}{\delta t}\|^2 + \frac{1}{2}\|\frac{e_u^1}{\delta t}\|^2 + \frac{1}{\delta t}\nu\|\nabla\tilde{e}_u^1\|^2 + \frac{1}{2}\|\nabla e_p^1\|^2 \\ = (\frac{\tilde{e}_u^1}{\delta t}, R_p^1) + (R_u^1, \frac{\tilde{e}_u^1}{\delta t}) + \frac{1}{2}\|R_p^1\|^2 - s(\boldsymbol{B}_0 \times \nabla \times e_b^1, \frac{\tilde{e}_u^1}{\delta t}).\end{aligned} \tag{3.71}$$

The terms on the right hand side can be esimated by

$$\begin{aligned}|(\frac{\tilde{e}_u^1}{\delta t}, R_p^1)| + |(R_u^1, \frac{\tilde{e}_u^1}{\delta t})| \lesssim \frac{1}{8}\|\frac{\tilde{e}_u^1}{\delta t}\|^2 + \|R_p^1\|^2 + \|R_u^1\|^2, \\ s|(\boldsymbol{B}_0 \times \nabla \times e_b^1, \frac{\tilde{e}_u^1}{\delta t})| \lesssim \frac{1}{8}\|\frac{\tilde{e}_u^1}{\delta t}\|^2 + \|\nabla \times e_b^1\|^2.\end{aligned} \tag{3.72}$$

From $e_u^0 = e_p^0 = 0$, thus (3.71) implies

$$\frac{1}{4}\|\frac{\tilde{e}_u^1}{\delta t}\|^2 + \frac{1}{2}\|\frac{e_u^1 - e_u^0}{\delta t}\|^2 + \frac{1}{\delta t}\nu\|\nabla\tilde{e}_u^1\|^2 + \frac{\delta t^2}{2}\|\frac{\nabla e_p^1 - \nabla e_p^0}{\delta t}\|^2 \lesssim \|R_p^1\|^2 + \|R_u^1\|^2 + \|\nabla \times e_b^1\|^2 \lesssim \delta t^2, \tag{3.73}$$

in which we use (3.68) and $\boldsymbol{u}_{tt} \in L^\infty(0,T;L^2), \boldsymbol{u}_t, p_t, \boldsymbol{B}_t \in L^\infty(0,T;H^1)$.

Finally, we obtain (3.65) by combining (3.68) and (3.73) together. ∎

**Lemma** 3.5. *Under assumption (A), there exists a constant $\widehat{C}$ such that, when $\delta t \le \widehat{C}$, the following estimate holds for $1 \le m \le [\frac{T}{\delta t}] - 1$,*

$$\|d_t e_u^{m+1}\|^2 + \|d_t e_b^{m+1}\|^2 + \delta t^2\|\nabla d_t e_p^{m+1}\|^2 + \delta t\sum_{n=1}^{m}(\nu\|\nabla d_t\tilde{e}_u^{n+1}\|^2 + \eta\|\nabla \times d_t e_b^{n+1}\|^2) \lesssim \delta t^2. \tag{3.74}$$

*Proof.* By taking the $L^2$ inner product of (3.61) with $d_t\tilde{e}_u^{n+1}$, we have

$$
\begin{aligned}
&\frac{1}{2\delta t}(\|d_t\tilde{e}_u^{n+1}-d_te_u^n\|^2+\|d_t\tilde{e}_u^{n+1}\|^2-\|d_te_u^n\|^2)+\nu\|\nabla d_t\tilde{e}_u^{n+1}\|^2\\
&\qquad+(\boldsymbol{u}(t_n)\nabla d_t\boldsymbol{u}(t_{n+1}),d_t\tilde{e}_u^{n+1})-(\boldsymbol{u}^n\nabla d_t\tilde{\boldsymbol{u}}^{n+1},d_t\tilde{e}_u^{n+1})\\
&\qquad+(d_t\boldsymbol{u}(t_n)\nabla\boldsymbol{u}(t_n),d_t\tilde{e}_u^{n+1})-(d_t\boldsymbol{u}^n\nabla\tilde{\boldsymbol{u}}^n,d_t\tilde{e}_u^{n+1})\\
&\qquad+s(\boldsymbol{B}(t_n)\times\nabla\times d_t\boldsymbol{B}(t_{n+1}),d_t\tilde{e}_u^{n+1})-s(\boldsymbol{B}^n\times\nabla\times d_t\boldsymbol{B}^{n+1},d_t\tilde{e}_u^{n+1})\\
&\qquad+s(d_t\boldsymbol{B}(t_n)\times\nabla\times\boldsymbol{B}(t_n),d_t\tilde{e}_u^{n+1})-s(d_t\boldsymbol{B}^n\times\nabla\times\boldsymbol{B}^n,d_t\tilde{e}_u^{n+1})\\
&\quad=(d_tR_u^{n+1},d_t\tilde{e}_u^{n+1})-(\nabla d_te_p^n,d_t\tilde{e}_u^{n+1}).
\end{aligned}\tag{3.75}
$$

From (3.62), we derive

$$
\frac{d_te_u^{n+1}}{\delta t}+\nabla d_te_p^{n+1}=\frac{d_t\tilde{e}_u^{n+1}}{\delta t}+\nabla d_te_p^n+d_tR_p^{n+1}.\tag{3.76}
$$

By taking the $L^2$ inner product of (3.76) with itself on both sides, we obtain

$$
\begin{aligned}
&\frac{1}{2\delta t}(\|d_te_u^{n+1}\|^2-\|d_t\tilde{e}_u^{n+1}\|^2)+\frac{\delta t}{2}(\|\nabla d_te_p^{n+1}\|^2-\|\nabla d_te_p^n\|^2)\\
&\quad=(d_t\tilde{e}_u^{n+1},d_tR_p^{n+1})+\delta t(\nabla d_te_p^n,d_tR_p^{n+1})\\
&\qquad+\frac{\delta t}{2}\|d_tR_p^{n+1}\|^2+(d_t\tilde{e}_u^{n+1},\nabla d_te_p^n).
\end{aligned}\tag{3.77}
$$

By taking the $L^2$ inner product of (3.63) with $d_te_b^{n+1}$, we have

$$
\begin{aligned}
&\frac{1}{2\delta t}(\|d_te_b^{n+1}-d_te_b^n\|^2+\|d_te_b^{n+1}\|^2-\|d_te_b^n\|^2)+\eta\|\nabla\times d_te_b^{n+1}\|^2\\
&\qquad+(d_t\boldsymbol{B}(t_n)\times\boldsymbol{u}(t_n),\nabla\times d_te_b^{n+1})-(d_t\boldsymbol{B}^n\times\boldsymbol{u}_\star^n,\nabla\times d_te_b^{n+1})\\
&\qquad+(\boldsymbol{B}(t_{n-1})\times d_t\boldsymbol{u}(t_n),\nabla\times d_te_b^{n+1})-(\boldsymbol{B}^{n-1}\times d_t\boldsymbol{u}_\star^n,\nabla\times d_te_b^{n+1})\\
&\qquad=(d_tR_b^{n+1},d_te_b^{n+1}).
\end{aligned}\tag{3.78}
$$

Combining (3.75), (3.77) and (3.78), we obtain

$$
\begin{aligned}
&\frac{1}{2\delta t}(\|d_t\tilde{e}_u^{n+1}-d_te_u^n\|^2+\|d_te_u^{n+1}\|^2-\|d_te_u^n\|^2)+\nu\|\nabla d_t\tilde{e}_u^{n+1}\|^2\\
&\qquad+\frac{\delta t}{2}(\|\nabla d_te_p^{n+1}\|^2-\|\nabla d_te_p^n\|^2)+\eta\|\nabla\times d_te_b^{n+1}\|^2\\
&\qquad+\frac{1}{2\delta t}(\|d_te_b^{n+1}-d_te_b^n\|^2+\|d_te_b^{n+1}\|^2-\|d_te_b^n\|^2)\\
&\qquad+(\boldsymbol{u}(t_n)\nabla d_t\boldsymbol{u}(t_{n+1}),d_t\tilde{e}_u^{n+1})-(\boldsymbol{u}^n\nabla d_t\tilde{\boldsymbol{u}}^{n+1},d_t\tilde{e}_u^{n+1}) &&(:\text{term I})\\
&\qquad+(d_t\boldsymbol{u}(t_n)\nabla\boldsymbol{u}(t_n),d_t\tilde{e}_u^{n+1})-(d_t\boldsymbol{u}^n\nabla\tilde{\boldsymbol{u}}^n,d_t\tilde{e}_u^{n+1}) &&(:\text{term II})\\
&\qquad+s(\boldsymbol{B}(t_n)\times\nabla\times d_t\boldsymbol{B}(t_{n+1}),d_t\tilde{e}_u^{n+1})-s(\boldsymbol{B}^n\times\nabla\times d_t\boldsymbol{B}^{n+1},d_t\tilde{e}_u^{n+1}) &&(:\text{term III})\\
&\qquad+s(d_t\boldsymbol{B}(t_n)\times\nabla\times\boldsymbol{B}(t_n),d_t\tilde{e}_u^{n+1})-s(d_t\boldsymbol{B}^n\times\nabla\times\boldsymbol{B}^n,d_t\tilde{e}_u^{n+1}) &&(:\text{term IV})\\
&\qquad+(d_t\boldsymbol{B}(t_n)\times\boldsymbol{u}(t_n),\nabla\times d_te_b^{n+1})-(d_t\boldsymbol{B}^n\times\boldsymbol{u}_\star^n,\nabla\times d_te_b^{n+1}) &&(:\text{term V})\\
&\qquad+(\boldsymbol{B}(t_{n-1})\times d_t\boldsymbol{u}(t_n),\nabla\times d_te_b^{n+1})-(\boldsymbol{B}^{n-1}\times d_t\boldsymbol{u}_\star^n,\nabla\times d_te_b^{n+1}) &&(:\text{term VI})\\
&\quad=(d_tR_u^{n+1},d_t\tilde{e}_u^{n+1})+(d_tR_b^{n+1},d_te_b^{n+1}) &&(:\text{term VII})\\
&\qquad+(d_t\tilde{e}_u^{n+1},d_tR_p^{n+1})+\delta t(\nabla d_te_p^n,d_tR_p^{n+1}) &&(:\text{term VIII})\\
&\qquad+\frac{\delta t}{2}\|d_tR_p^{n+1}\|^2.
\end{aligned}\tag{3.79}
$$

For term I, we estimate as

$$
\begin{aligned}
\text{(term I)} &\leq |(\boldsymbol{u}(t_n)\nabla d_t\boldsymbol{u}(t_{n+1}), d_t\tilde{e}_u^{n+1}) - (\boldsymbol{u}^n\nabla d_t\tilde{\boldsymbol{u}}^{n+1}, d_t\tilde{e}_u^{n+1})| \\
&= |(e_u^n\nabla d_t\boldsymbol{u}(t_{n+1}), d_t\tilde{e}_u^{n+1}) + (\boldsymbol{u}^n\nabla d_t\tilde{e}_u^{n+1}, d_t\tilde{e}_u^{n+1})| \\
&= |(e_u^n\nabla d_t\tilde{e}_u^{n+1}, d_t\boldsymbol{u}(t_{n+1}))| \\
&\lesssim \|\nabla e_u^n\|\|d_t\boldsymbol{u}(t_{n+1})\|_{L^4}\|\nabla d_t\tilde{e}_u^{n+1}\| \lesssim \frac{\nu}{12}\|\nabla d_t\tilde{e}_u^{n+1}\|^2 + \|\nabla e_u^n\|^2.
\end{aligned} \tag{3.80}
$$

where we use $\boldsymbol{u}_t \in L^\infty(0,T;L^4)$ and (3.17).

For term II, we estimate as

$$
\begin{aligned}
\text{(term II)} &\leq |(d_t\boldsymbol{u}(t_n)\nabla\boldsymbol{u}(t_n), d_t\tilde{e}_u^{n+1}) - (d_t\boldsymbol{u}^n\nabla\tilde{\boldsymbol{u}}^n, d_t\tilde{e}_u^{n+1})| \\
&= |(d_te_u^n\nabla\boldsymbol{u}(t_n), d_t\tilde{e}_u^{n+1}) + (d_t\boldsymbol{u}^n\nabla\tilde{e}_u^n, d_t\tilde{e}_u^{n+1})| \\
&= |(d_te_u^n\nabla\boldsymbol{u}(t_n), d_t\tilde{e}_u^{n+1}) - (d_te_u^n\nabla\tilde{e}_u^n, d_t\tilde{e}_u^{n+1}) + (d_t\boldsymbol{u}(t_n)\nabla\tilde{e}_u^n, d_t\tilde{e}_u^{n+1})| \\
&\lesssim \|d_te_u^n\|\|\nabla\boldsymbol{u}(t_n)\|_{L^4}\|\nabla d_t\tilde{e}_u^{n+1}\| + (\|d_te_u^n\|_{L^4} + \|d_t\boldsymbol{u}(t_n)\|_{L^4})\|\nabla\tilde{e}_u^n\|\|\nabla d_t\tilde{e}_u^{n+1}\| \\
&\lesssim \frac{\nu}{24}\|\nabla d_t\tilde{e}_u^{n+1}\|^2 + \|d_te_u^n\|^2 + \|d_te_u^n\|^{\frac{1}{2}}\|\nabla d_te_u^n\|^{\frac{1}{2}}\|\nabla\tilde{e}_u^n\|\|\nabla d_t\tilde{e}_u^{n+1}\| + \|d_t\boldsymbol{u}(t_n)\|_{L^4}^2\|\nabla\tilde{e}_u^n\|^2 \\
&\lesssim \frac{\nu}{12}\|\nabla d_t\tilde{e}_u^{n+1}\|^2 + \|d_te_u^n\|^2 + \|d_te_u^n\|\|\nabla d_te_u^n\|\|\nabla\tilde{e}_u^n\|^2 + \|d_t\boldsymbol{u}(t_n)\|_{L^4}^2\|\nabla\tilde{e}_u^n\|^2 \\
&\lesssim \frac{\nu}{12}\|\nabla d_t\tilde{e}_u^{n+1}\|^2 + \|d_te_u^n\|^2 + \|\nabla\tilde{e}_u^n\|^2 + \|d_te_u^n\|^2\|\nabla d_te_u^n\|^2\|\nabla\tilde{e}_u^n\|^2 + \|d_t\boldsymbol{u}(t_n)\|_{L^4}^2\|\nabla\tilde{e}_u^n\|^2 \\
&\lesssim \frac{\nu}{12}\|\nabla d_t\tilde{e}_u^{n+1}\|^2 + \|d_te_u^n\|^2 + \|\nabla\tilde{e}_u^n\|^2,
\end{aligned} \tag{3.81}
$$

where we actually use the fact that $\boldsymbol{u}_t \in L^\infty(0,T;L^4)$, and $\|\nabla d_te_u^n\|\|\nabla\tilde{e}_u^n\|$ is bounded since, from (3.59), we have

$$
\|\nabla d_te_u^n\|\|\nabla\tilde{e}_u^n\| \leq \frac{1}{\delta t}(\|\nabla e_u^n\| + \|\nabla e_u^{n-1}\|)\|\nabla\tilde{e}_u^n\| \lesssim \frac{1}{\delta t}(\|\nabla\tilde{e}_u^n\|^2 + \|\nabla\tilde{e}_u^{n-1}\|^2),
$$

that is bounded from Theorem 3.2.

For term III, we estimate as

$$
\begin{aligned}
\text{(term III)} &\leq |(\boldsymbol{B}(t_n)\times\nabla\times d_t\boldsymbol{B}(t_{n+1}), d_t\tilde{e}_u^{n+1}) - (\boldsymbol{B}^n\times\nabla\times d_t\boldsymbol{B}^{n+1}, d_t\tilde{e}_u^{n+1})| \\
&= |(e_b^n\times\nabla\times d_t\boldsymbol{B}(t_{n+1}), d_t\tilde{e}_u^{n+1}) + (\boldsymbol{B}^n\times\nabla\times d_te_b^{n+1}, d_t\tilde{e}_u^{n+1})| \\
&= |(\nabla\times(d_t\tilde{e}_u^{n+1}\times e_b^n), d_t\boldsymbol{B}(t_{n+1})) - (\nabla\times(\boldsymbol{B}^n\times d_t\tilde{e}_u^{n+1}), d_te_b^{n+1})| \\
&= |(e_b^n\nabla d_t\tilde{e}_u^{n+1} - d_t\tilde{e}_u^{n+1}\nabla e_b^n - e_b^n\nabla\cdot d_t\tilde{e}_u^{n+1}, d_t\boldsymbol{B}(t_{n+1})) \\
&\quad - (d_t\tilde{e}_u^{n+1}\nabla\boldsymbol{B}^n - \boldsymbol{B}^n\nabla d_t\tilde{e}_u^{n+1} + \boldsymbol{B}^n\nabla\cdot d_t\tilde{e}_u^{n+1}, d_te_b^{n+1})| \\
&\lesssim \|\nabla e_b^n\|\|\nabla d_t\tilde{e}_u^{n+1}\|\|d_t\boldsymbol{B}(t_{n+1})\|_{L^4} + \|\nabla d_t\tilde{e}_u^{n+1}\|\|\nabla e_b^n\|\|d_t\boldsymbol{B}(t_{n+1})\|_{L^4} \\
&\qquad + (\|\nabla d_t\tilde{e}_u^{n+1}\|\|\nabla\boldsymbol{B}^n\|_{L^4} + \|\boldsymbol{B}^n\|_{L^\infty}\|\nabla d_t\tilde{e}_u^{n+1}\|)\|d_te_b^{n+1}\| \\
&\lesssim \frac{\nu}{12}\|\nabla d_t\tilde{e}_u^{n+1}\|^2 + \|\nabla e_b^n\|^2\|d_t\boldsymbol{B}(t_{n+1})\|_{L^4}^2 + \|\nabla\boldsymbol{B}^n\|_{L^4}^2\|d_te_b^{n+1}\|^2 + \|\boldsymbol{B}^n\|_{L^\infty}^2\|d_te_b^{n+1}\|^2 \\
&\lesssim \frac{\nu}{12}\|\nabla d_t\tilde{e}_u^{n+1}\|^2 + \|\nabla\times e_b^n\|^2 + \|d_te_b^{n+1}\|^2,
\end{aligned} \tag{3.82}
$$

where we use $\boldsymbol{B}_t \in L^\infty(0,T;L^4)$ and $\|\nabla e_b^n\|^2 = \|\nabla\cdot e_b^n\|^2 + \|\nabla\times e_b^n\|^2$ and $\nabla\cdot e_b^n = 0$.

For term IV, we estimate as

$$
\begin{aligned}
(\text{term IV}) &\le |(d_t\boldsymbol{B}(t_n)\times\nabla\times\boldsymbol{B}(t_n), d_t\tilde{e}_u^{n+1}) - (d_t\boldsymbol{B}^n\times\nabla\times\boldsymbol{B}^n, d_t\tilde{e}_u^{n+1})| \\
&= |(d_te_b^n\times\nabla\times\boldsymbol{B}(t_n), d_t\tilde{e}_u^{n+1}) + (d_t\boldsymbol{B}^n\times\nabla\times e_b^n, d_t\tilde{e}_u^{n+1})| \\
&\lesssim \|d_te_b^n\|\|\boldsymbol{B}(t_n)\|_2\|\nabla d_t\tilde{e}_u^{n+1}\| + \|d_t\boldsymbol{B}^n\|_{L^4}\|\nabla\times e_b^n\|\|d_t\tilde{e}_u^{n+1}\|_{L^4} \\
&\lesssim \|d_te_b^n\|\|\nabla d_t\tilde{e}_u^{n+1}\| + (\|d_te_b^n\|_{L^4} + \|d_t\boldsymbol{B}(t_n)\|_{L^4})\|\nabla\times e_b^n\|\|\nabla d_t\tilde{e}_u^{n+1}\| \\
&\lesssim \frac{\nu}{12}\|\nabla d_t\tilde{e}_u^{n+1}\|^2 + \|d_te_b^n\|^2 + \|d_te_b^n\|_{L^4}^2\|\nabla\times e_b^n\|^2 + \|d_t\boldsymbol{B}(t_n)\|_{L^4}^2\|\nabla\times e_b^n\|^2 \\
&\lesssim \frac{\nu}{12}\|\nabla d_t\tilde{e}_u^{n+1}\|^2 + \|d_te_b^n\|^2 + \|d_te_b^n\|\|\nabla\times d_te_b^n\|\|\nabla\times e_b^n\|^2 + \|\nabla\times e_b^n\|^2 \\
&\lesssim \frac{\nu}{12}\|\nabla d_t\tilde{e}_u^{n+1}\|^2 + \|d_te_b^n\|^2 + \|d_te_b^n\|^2\|\nabla\times d_te_b^n\|^2\|\nabla\times e_b^n\|^2 + \|\nabla\times e_b^n\|^2 \\
&\lesssim \frac{\nu}{12}\|\nabla d_t\tilde{e}_u^{n+1}\|^2 + \|d_te_b^n\|^2 + \|\nabla\times e_b^n\|^2,
\end{aligned} \tag{3.83}
$$

where we use $\boldsymbol{B}_t\in L^\infty(0,T;L^4), \boldsymbol{B}\in L^\infty(0,T;H^2)$ and a fact that $\|\nabla\times d_te_b^n\|\|\nabla\times e_b^n\|$ is bounded since

$$
\|\nabla\times d_te_b^n\|\|\nabla\times e_b^n\| \le \frac{1}{\delta t}(\|\nabla\times e_b^n\| + \|\nabla\times e_b^{n-1}\|)\|\nabla\times e_b^n\| \lesssim \frac{1}{\delta t}(\|\nabla\times e_b^n\|^2 + \|\nabla\times e_b^{n-1}\|^2),
$$

that is bounded from Theorem 3.2.

For term V, we estimate as

$$
\begin{aligned}
(\text{term V}) &\le |(d_t\boldsymbol{B}(t_n)\times\boldsymbol{u}(t_n), \nabla\times d_te_b^{n+1}) - (d_t\boldsymbol{B}^n\times\boldsymbol{u}_\star^n, \nabla\times d_te_b^{n+1})| \\
&= |(d_te_b^n\times\boldsymbol{u}(t_n), \nabla\times d_te_b^{n+1}) + (d_t\boldsymbol{B}^n\times(\boldsymbol{u}(t_n)-\boldsymbol{u}_\star^n), \nabla\times d_te_b^{n+1})| \\
&= |(d_te_b^n\times\boldsymbol{u}(t_n), \nabla\times d_te_b^{n+1}) + ((d_t\boldsymbol{B}(t_n)-d_te_b^n)\times(\boldsymbol{u}(t_n)-\boldsymbol{u}_\star^n), \nabla\times d_te_b^{n+1})| \\
&\lesssim \|d_te_b^n\|\|\boldsymbol{u}(t_n)\|_{L^\infty}\|\nabla\times d_te_b^{n+1}\| + \|d_te_b^n\|_{L^4}\|\boldsymbol{u}(t_n)-\boldsymbol{u}_\star^n\|_{L^4}\|\nabla\times d_te_b^{n+1}\| \\
&\quad + \|d_t\boldsymbol{B}(t_n)\|_{L^4}\|\boldsymbol{u}(t_n)-\boldsymbol{u}_\star^n\|_{L^4}\|\nabla\times d_te_b^{n+1}\| \\
&\lesssim \|d_te_b^n\|\|\nabla\times d_te_b^{n+1}\| + \|d_te_b^n\|^{\frac12}\|\nabla\times d_te_b^n\|^{\frac12}\|\boldsymbol{u}(t_n)-\boldsymbol{u}_\star^n\|_1\|\nabla\times d_te_b^{n+1}\| \\
&\quad + \|\boldsymbol{u}(t_n)-\boldsymbol{u}_\star^n\|_1\|\nabla\times d_te_b^{n+1}\| \\
&\le \frac{\eta}{6}\|\nabla\times d_te_b^{n+1}\|^2 + \|d_te_b^n\|^2 + \|d_te_b^n\|\|\nabla\times d_te_b^n\|\|\boldsymbol{u}(t_n)-\boldsymbol{u}_\star^n\|_1^2 + \|\boldsymbol{u}(t_n)-\boldsymbol{u}_\star^n\|_1^2 \\
&\le \frac{\eta}{6}\|\nabla\times d_te_b^{n+1}\|^2 + \|d_te_b^n\|^2 + \|d_te_b^n\|^2\|\nabla\times d_te_b^n\|^2\|\boldsymbol{u}(t_n)-\boldsymbol{u}_\star^n\|_1^2 + \|\boldsymbol{u}(t_n)-\boldsymbol{u}_\star^n\|_1^2,
\end{aligned} \tag{3.84}
$$

where we use $\boldsymbol{u}\in L^\infty(0,T;H^2), \boldsymbol{B}_t\in L^\infty(0,T,H^1)$. In fact, by using (3.9) and (3.55), we get

$$
\begin{aligned}
\|\boldsymbol{u}(t_n)-\boldsymbol{u}_\star^n\| &\le \|e_u^n\| + s\delta t\|\boldsymbol{B}^n\times\nabla\times\boldsymbol{B}^{n+1}\| \\
&\lesssim \|e_u^n\| + \delta t\|\boldsymbol{B}^n\|_{L^\infty}\|\nabla\times\boldsymbol{B}^{n+1}\| \lesssim \delta t, \\
\|\nabla(\boldsymbol{u}(t_n)-\boldsymbol{u}_\star^n)\| &\le \|\nabla e_u^n\| + s\delta t\|\nabla(\boldsymbol{B}^n\times\nabla\times\boldsymbol{B}^{n+1})\| \\
&\lesssim \|\nabla e_u^n\| + \delta t(\|\boldsymbol{B}^n\|_2\|\boldsymbol{B}^{n+1}\|_2 + \|\boldsymbol{B}^n\|_{L^\infty}\|\boldsymbol{B}^{n+1}\|_2) \\
&\lesssim \|\nabla e_u^n\| + \delta t,
\end{aligned} \tag{3.85}
$$

which also implies

$$
\delta t\sum_{n=1}^{N}\|\boldsymbol{u}(t_n)-\boldsymbol{u}_\star^n\|_1^2 \lesssim \delta t^2. \tag{3.86}
$$

Therefore, from (3.86) and Theorem 3.2, we easily know that $\|\nabla\times d_te_b^n\|\|\boldsymbol{u}(t_n)-\boldsymbol{u}_\star^n\|_1 \lesssim \frac{1}{\delta t}(\|\nabla\times e_b^n\| + \|\nabla\times e_b^{n-1}\|)\|\boldsymbol{u}(t_n)-\boldsymbol{u}_\star^n\|_1$ is bounded by a constant, thus term V can be further estimated as

$$
(\text{term V}) \lesssim \frac{\eta}{6}\|\nabla\times d_te_b^{n+1}\|^2 + \|d_te_b^n\|^2 + \|\boldsymbol{u}(t_n)-\boldsymbol{u}_\star^n\|_1^2. \tag{3.87}
$$

For term VI, we estimate as

$$
\begin{aligned}
(\text{term VI}) &\le |(\boldsymbol{B}(t_{n-1}) \times d_t\boldsymbol{u}(t_n), \nabla\times d_te_b^{n+1}) - (\boldsymbol{B}^{n-1}\times d_t\boldsymbol{u}_\star^n, \nabla\times d_te_b^{n+1})| \\
&= |(e_b^{n-1}\times d_t\boldsymbol{u}(t_n), \nabla\times d_te_b^{n+1}) + (\boldsymbol{B}^{n-1}\times d_t(\boldsymbol{u}(t_n)-\boldsymbol{u}_\star^n), \nabla\times d_te_b^{n+1})| \\
&= |(e_b^{n-1}\times d_t\boldsymbol{u}(t_n), \nabla\times d_te_b^{n+1}) + \big(\boldsymbol{B}^{n-1}\times d_t(e_u^n + s\delta t\boldsymbol{B}^n\times\nabla\times\boldsymbol{B}^{n+1}), \nabla\times d_te_b^{n+1}\big)\,| \\
&\lesssim \|\nabla\times e_b^{n-1}\|\|d_t\boldsymbol{u}(t_n)\|_{L^4}\|\nabla\times d_te_b^{n+1}\| + \|\boldsymbol{B}^{n-1}\|_{L^\infty}\|d_te_u^n\|\|\nabla\times d_te_b^{n+1}\| \\
&\quad + \delta t|(\boldsymbol{B}^{n-1}\times d_t(\boldsymbol{B}^n\times\nabla\times\boldsymbol{B}^{n+1}), \nabla\times d_te_b^{n+1})| \\
&\lesssim \|\nabla\times e_b^{n-1}\|\|\nabla\times d_te_b^{n+1}\| + \|d_te_u^n\|\|\nabla\times d_te_b^{n+1}\| + \delta t|(\boldsymbol{B}^{n-1}\times d_t(\boldsymbol{B}^n\times\nabla\times\boldsymbol{B}^{n+1}), \nabla\times d_te_b^{n+1})| \\
&\lesssim \frac{\eta}{12}\|\nabla\times d_te_b^{n+1}\|^2 + \|\nabla\times e_b^{n-1}\|^2 + \|d_te_u^n\|^2 + \delta t|(\boldsymbol{B}^{n-1}\times d_t(\boldsymbol{B}^n\times\nabla\times\boldsymbol{B}^{n+1}), \nabla\times d_te_b^{n+1})|,
\end{aligned}
\tag{3.88}
$$

where we use $\boldsymbol{u}_t \in L^\infty(0,T;H^1(\Omega))$. The last term on the right hand side of (3.88) can be estimated as

$$
\begin{aligned}
&\delta t|(\boldsymbol{B}^{n-1}\times d_t(\boldsymbol{B}^n\times\nabla\times\boldsymbol{B}^{n+1}), \nabla\times d_te_b^{n+1})| \\
&\quad\lesssim \|\boldsymbol{B}^{n-1}\|_{L^\infty}\|\boldsymbol{B}^n\times\nabla\times\boldsymbol{B}^{n+1} - \boldsymbol{B}^{n-1}\times\nabla\times\boldsymbol{B}^n)\|\|\nabla\times d_te_b^{n+1}\| \\
&\quad\lesssim \|\boldsymbol{B}^n\times\nabla\times(\boldsymbol{B}^{n+1}-\boldsymbol{B}^n) + (\boldsymbol{B}^n-\boldsymbol{B}^{n-1})\times\nabla\times\boldsymbol{B}^n\|\|\nabla\times d_te_b^{n+1}\| \\
&\quad\lesssim \big(\|\boldsymbol{B}^n\times\nabla\times(\boldsymbol{B}^{n+1}-\boldsymbol{B}^n)\| + \|(\boldsymbol{B}^n-\boldsymbol{B}^{n-1})\times\nabla\times\boldsymbol{B}^n\|\big)\,\|\nabla\times d_te_b^{n+1}\| \\
&\quad\lesssim \|\boldsymbol{B}^n\times\nabla\times(e_b^n - e_b^{n+1} + \boldsymbol{B}(t_{n+1}) - \boldsymbol{B}(t_n))\|\|\nabla\times d_te_b^{n+1}\| \\
&\qquad + \|(e_b^{n-1} - e_b^n + \boldsymbol{B}(t_n) - \boldsymbol{B}(t_{n-1}))\times\nabla\times\boldsymbol{B}^n\|\|\nabla\times d_te_b^{n+1}\| \\
&\quad\lesssim \big(\|\nabla\times(e_b^n - e_b^{n+1})\| + \|\nabla\times(\boldsymbol{B}(t_{n+1}) - \boldsymbol{B}(t_n))\|\big)\,\|\nabla\times d_te_b^{n+1}\| \\
&\qquad + \big(\|e_b^{n-1} - e_b^n\|_{L^4}\|\boldsymbol{B}^n\|_2 + \|\boldsymbol{B}(t_n) - \boldsymbol{B}(t_{n-1})\|_1\|\boldsymbol{B}^n\|_2\big)\,\|\nabla\times d_te_b^{n+1}\| \\
&\quad\lesssim \frac{\eta}{12}\|\nabla\times d_te_b^{n+1}\|^2 + \|\nabla\times(e_b^n - e_b^{n+1})\|^2 + \|\boldsymbol{B}(t_{n+1}) - \boldsymbol{B}(t_n)\|_1^2 \\
&\qquad + \|\nabla\times(e_b^{n-1} - e_b^n)\|^2 + \|\boldsymbol{B}(t_n) - \boldsymbol{B}(t_{n-1})\|_1^2,
\end{aligned}
\tag{3.89}
$$

where we use Remark 3.1 and Lemma 3.2. Therefore, for term VI, we obtain

$$
\begin{aligned}
(\text{term VI}) &\lesssim \frac{\eta}{6}\|\nabla\times d_te_b^{n+1}\|^2 + \|\nabla\times e_b^{n-1}\|^2 + \|d_te_u^n\|^2 + \|\nabla\times(e_b^n - e_b^{n+1})\|^2 + \|\nabla\times(e_b^{n-1} - e_b^n)\|^2 \\
&\quad + \|\boldsymbol{B}(t_{n+1}) - \boldsymbol{B}(t_n)\|_1^2 + \|\boldsymbol{B}(t_n) - \boldsymbol{B}(t_{n-1})\|_1^2.
\end{aligned}
\tag{3.90}
$$

For term VII, we have

$$
\begin{aligned}
(\text{term VII}) &\le |(d_tR_u^{n+1}, d_t\tilde{e}_u^{n+1}) + (d_tR_b^{n+1}, d_te_b^{n+1})| \\
&\lesssim \frac{\nu}{12}\|\nabla d_t\tilde{e}_u^{n+1}\|^2 + \frac{\eta}{6}\|\nabla\times d_te_b^{n+1}\|^2 + \|d_tR_u^{n+1}\|^2 + \|d_tR_b^{n+1}\|^2.
\end{aligned}
\tag{3.91}
$$

For term VIII, we have

$$
\begin{aligned}
(\text{term VIII}) &\le |(d_t\tilde{e}_u^{n+1}, d_tR_p^{n+1}) - \delta t(\nabla d_te_p^n, d_tR_p^{n+1})| \\
&\lesssim \frac{\nu}{12}\|\nabla d_t\tilde{e}_u^{n+1}\|^2 + \|d_tR_p^{n+1}\|^2 + \delta t^2\|\nabla d_te_p^n\|^2.
\end{aligned}
\tag{3.92}
$$

By combining (3.79), (3.80), (3.81), (3.82), (3.83), (3.84), (3.90), (3.91) and (3.92), we have

$$
\begin{aligned}
&\|d_te_u^{n+1}\|^2 - \|d_te_u^n\|^2 + \delta t\nu\|\nabla d_t\tilde{e}_u^{n+1}\|^2 + \delta t^2(\|\nabla d_te_p^{n+1}\|^2 - \|\nabla d_te_p^n\|^2) \\
&\qquad + \|d_te_b^{n+1}\|^2 - \|d_te_b^n\|^2 + \delta t\eta\|\nabla\times d_te_b^{n+1}\|^2 \\
&\qquad \lesssim \delta t(\|\nabla e_u^n\|^2 + \|\nabla\tilde{e}_u^n\|^2 + \|\nabla\times e_b^n\|^2 + \|\nabla\times e_b^{n-1}\|^2 + \|\nabla\times e_b^{n+1}\|^2) \\
&\qquad + \delta t(\|d_te_u^n\|^2 + \|d_te_b^n\|^2 + \|d_te_b^{n+1}\|^2 + \delta t^2\|\nabla d_te_p^n\|^2) \\
&\qquad + \delta t(\|d_tR_p^{n+1}\|^2 + \|d_tR_u^{n+1}\|^2 + \|d_tR_b^{n+1}\|^2) \\
&\qquad + \delta t\big(\|\boldsymbol{B}(t_{n+1}) - \boldsymbol{B}(t_n)\|_1^2 + \|\boldsymbol{B}(t_n) - \boldsymbol{B}(t_{n-1})\|_1^2\big) \\
&\qquad + \delta t\|\boldsymbol{u}(t_n) - \boldsymbol{u}_\star^n\|_1^2.
\end{aligned}
\tag{3.93}
$$

By taking the summation of (3.93) from $n = 1$ to $m$, we obtain

$$
\begin{aligned}
&\|d_t e_u^{m+1}\|^2 + \|d_t e_b^{m+1}\|^2 + \delta t^2 \|\nabla d_t e_p^{m+1}\|^2 + \delta t \sum_{n=1}^{m} (\nu \|\nabla d_t \tilde{e}_u^{n+1}\|^2 + \eta \|\nabla \times d_t e_b^{n+1}\|^2) \\
&\lesssim \|d_t e_u^1\|^2 + \|d_t e_b^1\|^2 + \delta t^2 \|\nabla d_t e_p^1\|^2 + \delta t \sum_{n=1}^{m} (\|d_t e_u^n\|^2 + \|d_t e_b^n\|^2 + \|d_t e_b^{n+1}\|^2 + \delta t^2 \|\nabla d_t e_p^n\|^2) \\
&+ \delta t \sum_{n=1}^{m} (\|d_t R_p^{n+1}\|^2 + \|d_t R_u^{n+1}\|^2 + \|d_t R_b^{n+1}\|^2) + \delta t^2 \\
&\lesssim \delta t \sum_{n=1}^{m} (\|d_t e_u^n\|^2 + \|d_t e_b^n\|^2 + \|d_t e_b^{n+1}\|^2 + \delta t^2 \|\nabla d_t e_p^n\|^2) + \delta t^2,
\end{aligned} \tag{3.94}
$$

where we use Lemma 3.3, Lemma 3.4, (3.60), (3.86) and the assumption (A).

From the Gronwall's inequality, there exists a constant $\widehat{C}$, such that

$$
\|d_t e_u^{m+1}\|^2 + \|d_t e_b^{m+1}\|^2 + \delta t^2 \|\nabla d_t e_p^{m+1}\|^2 + \delta t \sum_{n=1}^{m} (\nu \|\nabla d_t \tilde{e}_u^{n+1}\|^2 + \eta \|\nabla \times d_t e_b^{n+1}\|^2) \lesssim \delta t^2,
$$

holds for $\delta t \le \widehat{C}$ and $1 \le m \le [\frac{T}{\delta t}] - 1$, that concludes (3.74). □

**Theorem 3.3.** *Under the assumptions of Lemma 3.5, we have*

$$
\|\nabla e_u^n\| + \|\nabla \tilde{e}_u^n\| + \|\nabla \times e_b^n\| + \|e_p^n\| \lesssim \delta t, \ 1 \le n \le [\frac{T}{\delta t}]. \tag{3.95}
$$

*Proof.* By taking the summation of (3.29) and (3.30), we obtain the error equation as

$$
\begin{aligned}
&\frac{e_u^{n+1} - e_u^n}{\delta t} - \nu \Delta \tilde{e}_u^{n+1} + (\boldsymbol{u}(t_n) \cdot \nabla) \boldsymbol{u}(t_{n+1}) - (\boldsymbol{u}^n \cdot \nabla) \tilde{\boldsymbol{u}}^{n+1} + \nabla e_p^{n+1} \\
&+ s \boldsymbol{B}(t_n) \times \nabla \times \boldsymbol{B}(t_{n+1}) - s \boldsymbol{B}^n \times \nabla \times \boldsymbol{B}^{n+1} = R_u^{n+1} + R_p^{n+1}.
\end{aligned} \tag{3.96}
$$

By taking $L^2$ inner product of (3.96) with any test function $\boldsymbol{v}$ $(\in H_0^1(\Omega)^2)$, from assumption (A), we obtain

$$
\begin{aligned}
(e_p^{n+1}, \nabla \cdot \boldsymbol{v}) = &\left( \frac{e_u^{n+1} - e_u^n}{\delta t}, \boldsymbol{v} \right) + \nu (\nabla \tilde{e}_u^{n+1}, \nabla \boldsymbol{v}) + \left( (e_u^n \cdot \nabla) \boldsymbol{u}(t_{n+1}) + (\boldsymbol{u}^n \cdot \nabla) \tilde{e}_u^{n+1}, \boldsymbol{v} \right) \\
&+ s \left( e_b^n \times \nabla \times \boldsymbol{B}(t_{n+1}) + \boldsymbol{B}^n \times \nabla \times e_b^{n+1}, \boldsymbol{v} \right) - (R_u^{n+1} + R_p^{n+1}, \boldsymbol{v}) \\
&\lesssim \|d_t e_u^{n+1}\| \|\nabla \boldsymbol{v}\| + \|\nabla \tilde{e}_u^{n+1}\| \|\nabla \boldsymbol{v}\| + \|e_u^n\| \|\nabla \boldsymbol{v}\| + \|\nabla \tilde{e}_u^{n+1}\| \|\nabla \boldsymbol{v}\| \\
&+ \|e_b^n\| \|\nabla \boldsymbol{v}\| + \|\nabla \times e_b^{n+1}\| \|\nabla \boldsymbol{v}\| + \|R_u^{n+1}\| \|\nabla \boldsymbol{v}\| + \|R_p^{n+1}\| \|\nabla \boldsymbol{v}\|.
\end{aligned} \tag{3.97}
$$

Using the inf-sup condition

$$
\beta \|q\| \le \sup_{\boldsymbol{w} \in H_0^1(\Omega)} \frac{(\nabla \cdot \boldsymbol{w}, q)}{\|\nabla \boldsymbol{w}\|}, \quad \forall q \in L_0^2(\Omega), \tag{3.98}
$$

we obtain

$$
\begin{aligned}
\beta \|e_p^{n+1}\| &\lesssim \|d_t e_u^{n+1}\| + \|\nabla \tilde{e}_u^{n+1}\| + \|e_u^n\| + \|e_b^n\| + \|\nabla \times e_b^{n+1}\| + \|R_u^{n+1}\| + \|R_p^{n+1}\| \\
&\lesssim \delta t + \|\nabla \tilde{e}_u^{n+1}\| + \|\nabla \times e_b^{n+1}\|.
\end{aligned} \tag{3.99}
$$

From Lemma 3.5 and Hölder's inequality, for $1 \le m \le [\frac{T}{\delta t}] - 1$, we have

$$
\begin{aligned}
\|\nabla \tilde{e}_u^{m+1}\| - \|\nabla \tilde{e}_u^1\| &= \sum_{n=1}^{m} (\|\nabla \tilde{e}_u^{n+1}\| - \|\nabla \tilde{e}_u^n\|) \le \sum_{n=1}^{m} \|\nabla \tilde{e}_u^{n+1} - \nabla \tilde{e}_u^n\| \\
&\le (\delta t^2 \sum_{n=1}^{m} \|\frac{\nabla \tilde{e}_u^{n+1} - \nabla \tilde{e}_u^n}{\delta t}\|^2)^{\frac{1}{2}} (\sum_{n=1}^{m} 1^2)^{\frac{1}{2}} \lesssim (\delta t^3)^{\frac{1}{2}} (\frac{T}{\delta t})^{\frac{1}{2}} \lesssim \delta t.
\end{aligned} \tag{3.100}
$$

Thus, from (3.73), we obtain

$$\|\nabla \tilde{e}_u^{m+1}\| \lesssim \delta t + \|\nabla \tilde{e}_u^1\| \lesssim \delta t, 0 \le m \le [\frac{T}{\delta t}] - 1. \tag{3.101}$$

Due to the $H^1$ stability of $P_H$ [36], we have

$$\|\nabla e_u^{m+1}\| \le \|\nabla \tilde{e}_u^{m+1}\| \lesssim \delta t. \tag{3.102}$$

Similarly, we can also deduce

$$\|\nabla \times e_b^{m+1}\| \lesssim \delta t. \tag{3.103}$$

Finally, from (3.99), (3.101), (3.102), (3.103), Lemma 3.5 and Theorem 3.2, we have (3.95). ◻

**4. Numerical examples.** We now implement various 2D numerical experiments to validate the stability and accuracy of the scheme. We use the inf-sup stable P2/P1 element [5] for the velocity and pressure, and linear element for the magnetic field.

The fully discretized finite element scheme of (3.1)-(3.6) reads as

*Step I.* Find $\boldsymbol{B}_h^{n+1}$ such that

$$\left(\frac{\boldsymbol{B}_h^{n+1} - \boldsymbol{B}_h^n}{\delta t}, \boldsymbol{C}_h\right) + \eta(\nabla \times \boldsymbol{B}_h^{n+1}, \nabla \times \boldsymbol{C}_h) + \eta(\nabla \cdot \boldsymbol{B}_h^{n+1}, \nabla \cdot \boldsymbol{C}_h) + (\boldsymbol{B}_h^n \times \boldsymbol{u}_h^n, \nabla \times \boldsymbol{C}_h)$$
$$+ \delta t s(\boldsymbol{B}_h^n \times \nabla \times \boldsymbol{B}_h^{n+1}, \boldsymbol{B}_h^n \times \nabla \times \boldsymbol{C}_h) = 0.$$

*Step II.* Find $\tilde{\boldsymbol{u}}_h^{n+1}$ such that

$$\left(\frac{\tilde{\boldsymbol{u}}_h^{n+1} - \boldsymbol{u}_h^n}{\delta t}, \boldsymbol{v}_h\right) + \nu(\nabla \tilde{\boldsymbol{u}}_h^{n+1}, \nabla \boldsymbol{v}_h) + (\boldsymbol{u}_h^n \nabla \tilde{\boldsymbol{u}}_h^{n+1}, \boldsymbol{v}_h) - (p_h^n, \nabla \cdot \boldsymbol{v}_h) + s(\boldsymbol{B}_h^n \times \nabla \times \boldsymbol{B}_h^{n+1}, \boldsymbol{v}_h) = 0.$$

*Step III.* Find $p_h^{n+1}$ from

$$(\nabla p_h^{n+1}, \nabla q_h) = -\frac{1}{\delta t}(\nabla \cdot \tilde{\boldsymbol{u}}_h^{n+1}, q_h) + (\nabla p_h^n, \nabla q_h).$$

*Step IV.* Update $\boldsymbol{u}_h^{n+1}$ from

$$(\boldsymbol{u}_h^{n+1}, \boldsymbol{v}_h) = (\tilde{\boldsymbol{u}}_h^{n+1}, \boldsymbol{v}_h) + \delta t(p_h^{n+1}, \nabla \cdot \boldsymbol{v}_h) - \delta t(p_h^n, \nabla \cdot \boldsymbol{v}_h).$$

**Remark 4.1.** *It is remarkable that, in the fully discrete case, the divergence free condition for the magnetic field is not valid any more. Therefore a common practice is to add a penalty term* $(\nabla \cdot \boldsymbol{B}_h^{n+1}, \nabla \cdot \boldsymbol{C}_h)$ *in the momentum equation (cf. see [14]) in order to ensure the coercivity of the magnetic equation and approximate divergence free of magnetic field in the fully discrete level.*

**4.1. Accuracy test.** We first perform numerical simulations to test the convergence rates of the proposed scheme. The computational domain is $\Omega = [0,1] \times [0,1]$. We assume the following functions

$$\begin{cases} \boldsymbol{u} = (y^5 + t^2, x^5 + t^2), \\ p = 10(2x-1)(2y-1)(1+t^2), \\ \boldsymbol{B} = (\sin(y) + t^2, \sin(x) + t^2) \end{cases} \tag{4.1}$$

to be the exact solution, and impose some suitable force fields such that the given solution can satisfy the system. The physical parameters are $\nu = \eta = s = 1$. In Table 1, we list the numerical errors between the numerical solution and the exact solution at $T = 1$ with different time step sizes. Since the relation between the spatial grid size $h$ and temporal step size $\delta t$ to be $\delta t = h^2$, we observe the second order accuracy asymptotically for $\|e_u\|_{L^2}, \|e_u\|_{H^1}, \|e_p\|_{L^2}, \|e_b\|_{L^2}$ and the first order accuracy for $\|e_b\|_{H^1}$, as predicted theoretically.

| $\delta t$ | $h$ | $\|e_u\|_{L^2}$ | Order | $\|e_u\|_{H^1}$ | Order | $\|e_p\|_{L^2}$ | Order | $\|e_b\|_{L^2}$ | Order | $\|e_b\|_{H^1}$ | Order |
|---|---|---|---|---|---|---|---|---|---|---|---|
| 1/16 | 1/4 | 4.03e-2 | – | 5.43e-1 | – | 1.53 | – | 9.67e-3 | – | 7.17e-2 | – |
| 1/64 | 1/8 | 3.56e-3 | 3.50 | 8.22e-2 | 2.72 | 2.03e-1 | 2.91 | 3.50e-3 | 1.47 | 3.19e-2 | 1.16 |
| 1/256 | 1/16 | 2.47e-4 | 3.85 | 1.01e-2 | 3.02 | 2.63e-2 | 2.94 | 9.50e-4 | 1.88 | 1.41e-2 | 1.17 |
| 1/1024 | 1/32 | 2.25e-5 | 3.45 | 1.60e-3 | 2.65 | 5.30e-3 | 2.31 | 2.42e-4 | 1.97 | 6.77e-3 | 1.06 |
| 1/4096 | 1/64 | 4.30e-6 | 2.39 | 3.55e-4 | 2.17 | 1.30e-3 | 2.02 | 6.09e-05 | 1.99 | 3.34e-3 | 1.01 |
| 1/16384 | 1/128 | 1.05e-6 | 2.02 | 8.73e-5 | 2.02 | 3.25e-4 | 2.00 | 1.52e-5 | 2.00 | 1.66e-3 | 1.00 |

TABLE 1

*The numerical errors and convergence order for $\|e_u\|_{L^2}, \|e_u\|_{H^1}, \|e_p\|_{L^2}, \|e_b\|_{L^2}, \|e_b\|_{H^1}$ at $t = 1$ that are computed using various temporal and spatial resolutions with the exact solutions of (4.1). The physical parameters are $\nu = \eta = s = 1$, and the time step and the grid size satisfy $\delta t = h^2$.*

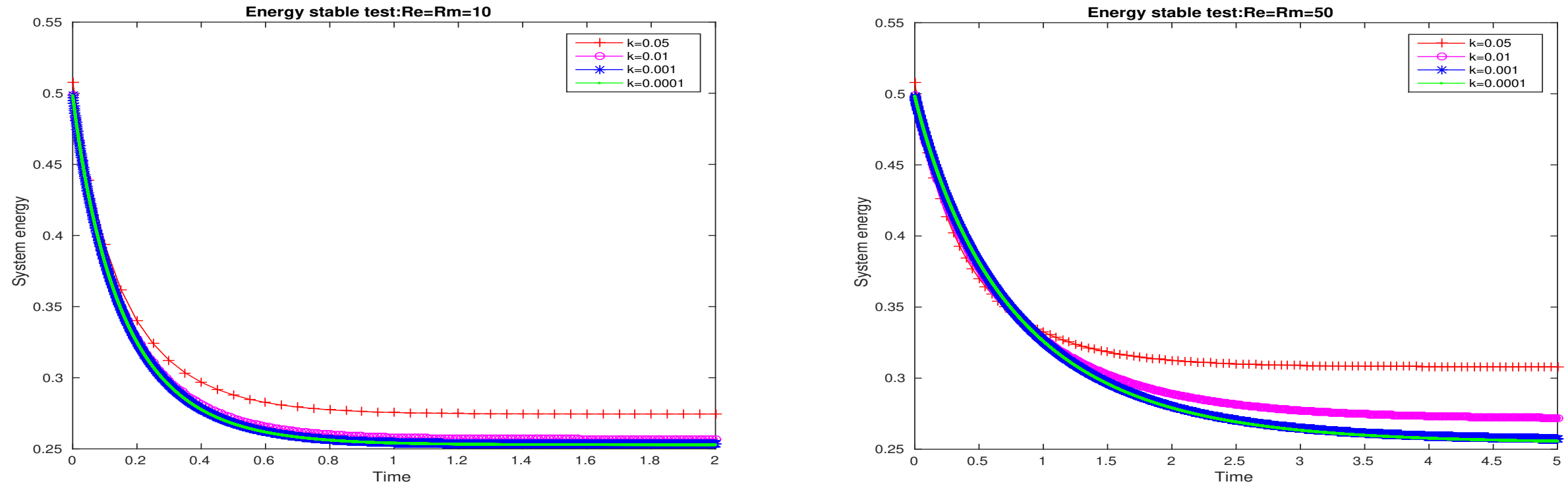

(a) $Re = Rm = 10$. (b) $Re = Rm = 50$.

FIG. 1. *Time evolution of the free energy functional till $T = 5$ for four different time steps and two sets of order parameters (a) $R_e = R_m = 10$ and (b) $R_e = R_m = 50$. The energy curves show the decays for all time steps $\delta t = 0.05, 0.01, 0.001, 0.0001$, which confirms that our algorithm is unconditionally stable.*

**4.2. Stability test.** We show the evolution of the total free energy in this example. We set the computed domain to be $\Omega = [0, 1]^2$, and the initial conditions for $\boldsymbol{u}, p, \boldsymbol{B}$ are

$$(4.2)\qquad \begin{cases} \boldsymbol{u}^0 = (x^2(x-1)^2 y(y-1)(2y-1),\ -y^2(y-1)^2 x(x-1)(2x-1)), p^0 = 0, \\ \boldsymbol{B}^0 = (\sin(\pi x)\cos(\pi y),\ -sin(\pi y)\cos(\pi x)). \end{cases}$$

We test the energy stability over matching time of the proposed scheme under variant physical parameters of $R_e = R_m = 10$ and 50. The coupling parameter is fixed as $s = 1$, and mesh size is $h = 1/64$. In Fig. 1, we present the time evolution of the total free energy for four different time steps of $k = 0.05, 0.01, 0.001, 0.0001$ until $T = 5$. We observe that all four energy curves show decays monotonically for all time step sizes, which numerically confirms that our algorithm is unconditionally energy stable.

**4.3. Hydromagnetic Kelvin-Helmholtz instability.** The Kelvin-Helmholtz (K-H) instability in sheared flow configurations is an efficient mechanism to initiate mixing of fluids, transport of momentum and energy, and the development of turbulence. Such a problem is of interest in investigating a variety of space, astrophysical, and geophysical situations involving sheared plasma flows. Configurations where it is relevant include the interface between the solar wind and the magnetosphere, coronal streamers moving through the solar wind, etc. Since most astrophysical environments are electrically conducting and relevant fluids are likely to be magnetized, it is thus of prime importance to understand the role of magnetic fields in the K-H instability. About the theoretical and numerical study of Hydromagnetic K-H instability, we refer to [28, 11, 9, 25, 8, 16, 4] and the references therein.

We revisit the occurring of the K-H instability in a single shear flow configuration that is embedded in a uniform flow-aligned magnetic field. The simulation is performed in the computed domain of $[0, 2] \times [0, 1]$. The initial velocity field is $\boldsymbol{u}_0 = (1.5, 0)$ in the top half domain, and $\boldsymbol{u}_0 = (-1.5, 0)$ in the bottom half domain. The sheared initial magnetic field is $\boldsymbol{B}_0 = (\tanh(y/\epsilon), 0)$ where $\epsilon = 0.07957747154595$ (cf. [9]). The velocity $\boldsymbol{u}$, magnetic field $\boldsymbol{B}$ and pressure $p$ are periodic boundary conditions on left and right boundaries. On the top and bottom boundary, the second component $v$ of the velocity field $\boldsymbol{u} = (u, v)$ is imposed. The boundary conditions for $\boldsymbol{B}$ are $\boldsymbol{B} \times \boldsymbol{n} = \boldsymbol{B}_0 \times \boldsymbol{n}$ for the top boundary and $-\boldsymbol{B}_0 \times \boldsymbol{n}$ for the bottom. The order parameters are $R_e = R_m = 1000$, $s = 0.2$. We use the time step $\delta t = 0.01$ and grid size $h = \frac{1}{40}$.

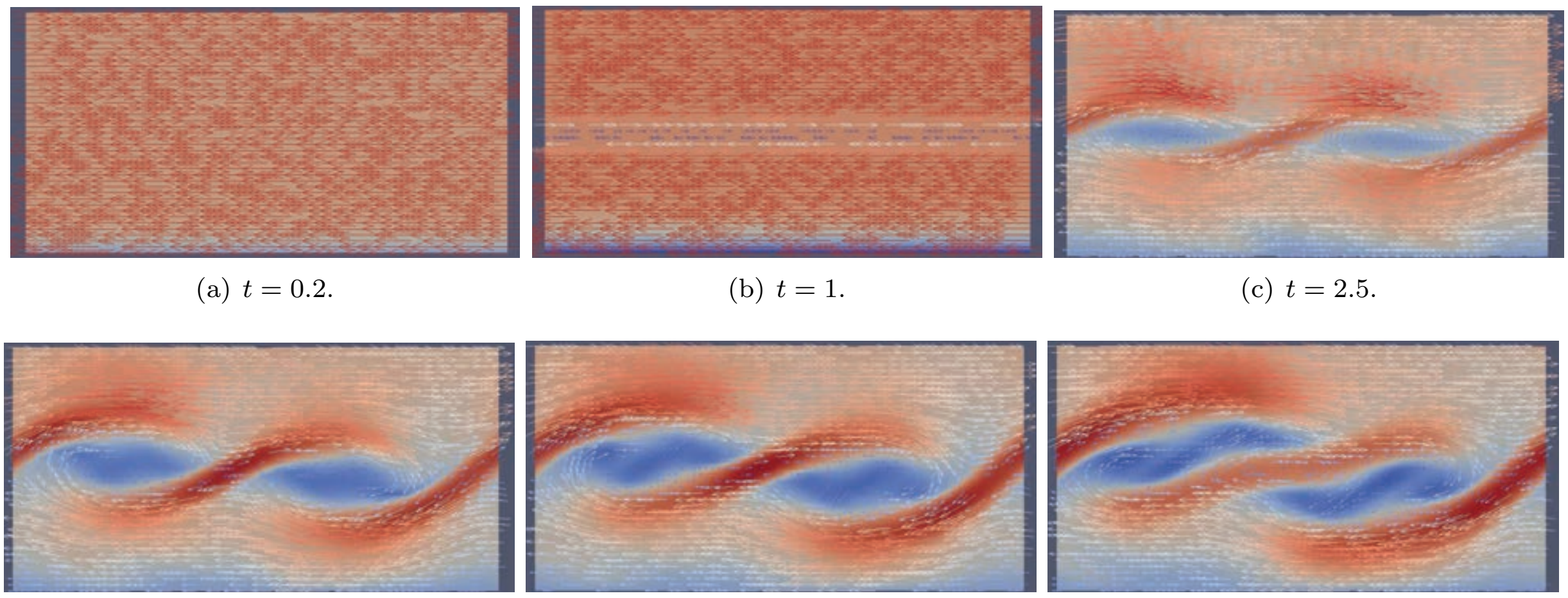

(a) $t = 0.2$. (b) $t = 1$. (c) $t = 2.5$.

(d) $t = 3$. (e) $t = 3.2$. (f) $t = 3.5$.

FIG. 2. *The dynamical behaviors of the magnetic field together with the velocity field that shows the hydromagnetic K-H instability. Snapshots of the numerical approximation are taken at* $t = 0.2, 1, 2.5, 3, 3.2, 3.5$.

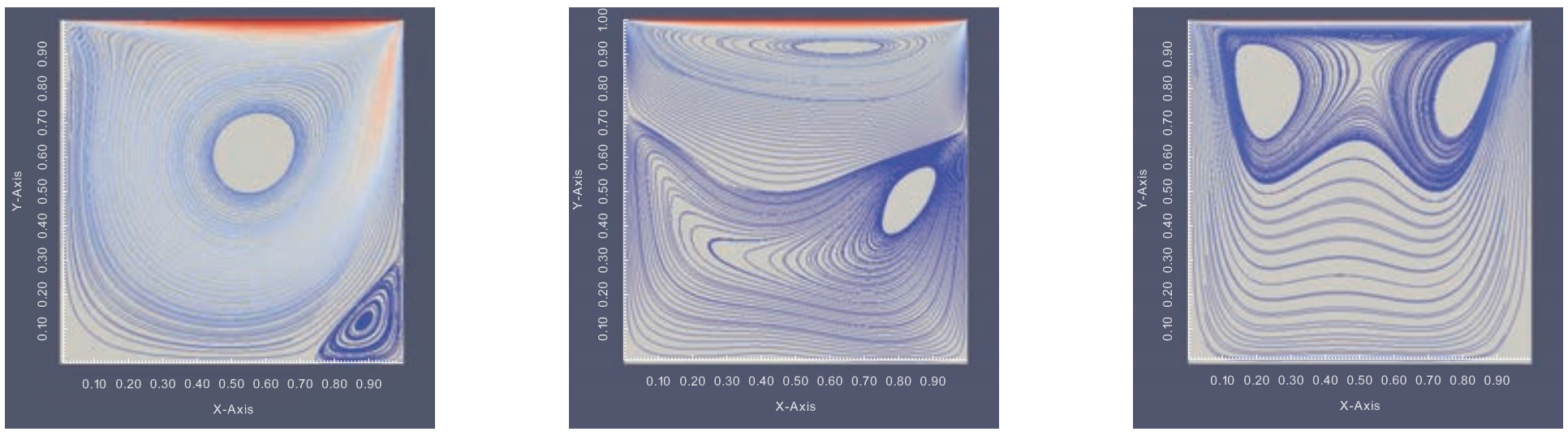

(a) No magnetic field. (b) H-Magnetic. (c) V-Magnetic.

FIG. 3. *The steady state of the streamlines for the velocity field for three cases, where (a) No magnetic field, (b) H-Magnetic, and (c) V-Magnetic.*

In Fig. 2, we show snapshots of the magnitude of $\boldsymbol{B}_1$ that is the first component of $\boldsymbol{B} = (\boldsymbol{B}_1, \boldsymbol{B}_2)$ together with the velocity field $\boldsymbol{u}$ at $t = 0.2, 1, 2.5, 3, 3.2, 3.5$. When time evolves, we can observe the vortexes start to form around $t = 1$. After $t = 2.5$, the profiles of vortexes and the magnetic field show the typical structure of K-H instability, and soon it deforms and rotates along with the flow. The obtained numerical results coincide well with the numerical/experimental results discussed in [28, 16, 8, 4], qualitatively.

**4.4. Lid driven cavity.** Finally we perform the lid driven cavity flow simulation of a conducting fluid. The computed domain is $(x, y) \in \Omega = [0, 1] \times [0, 1]$. The boundary condition of the magnetic field is either horizontal case $\boldsymbol{B}|_{\partial\Omega} = (1, 0)$ (denoted by H-magnetic) or vertical case $\boldsymbol{B}|_{\partial\Omega} = (0, 1)$ (denoted by V-magnetic). The no-slip boundary conditions are imposed on the bottom, left, and right sides of the cavity and the lid moves from left to right with the constant speed $\boldsymbol{u}|_{y=1} = (1, 0)$.

First, we set $R_e = 400, R_m = 40$ and Hartmann number $H_a(:= \sqrt{sR_eR_m}) = 100$. In Fig. 3, we show the steady state for three cases: no magnetic field, H-magnetic and V-magnetic, respectively. For the case of no magnetic field, this situation is actually a pure hydrodynamic problem, in which we observe the cavity is dominated by one large primary eddy in the center, together with a smaller secondary eddy in the right corner. For the H-magnetic case, two large eddies appear in the upper and lower part. For the V-magnetic case, two large eddies are located in the upper left and upper right. Similar features had been observed in [2] as well.

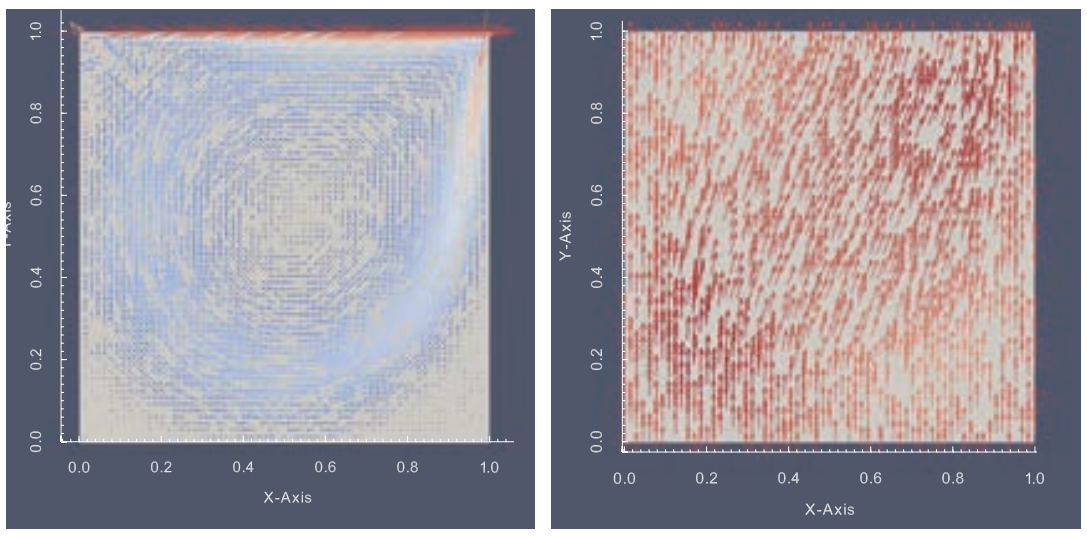

(a) $R_e = 1000, R_m = 10, H_a = 10$.

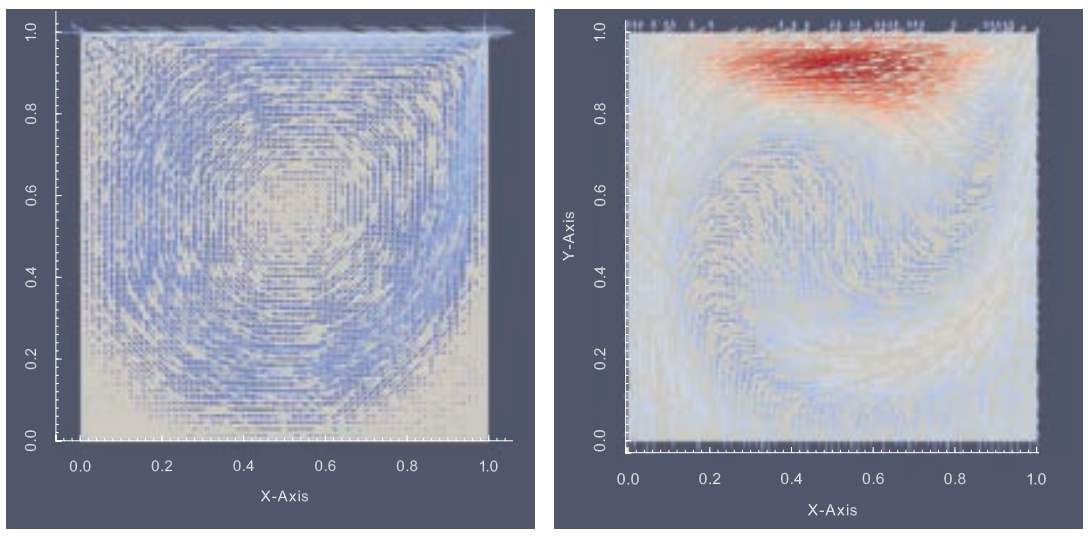

(b) $R_e = 1000, R_m = 100, H_a = 10$.

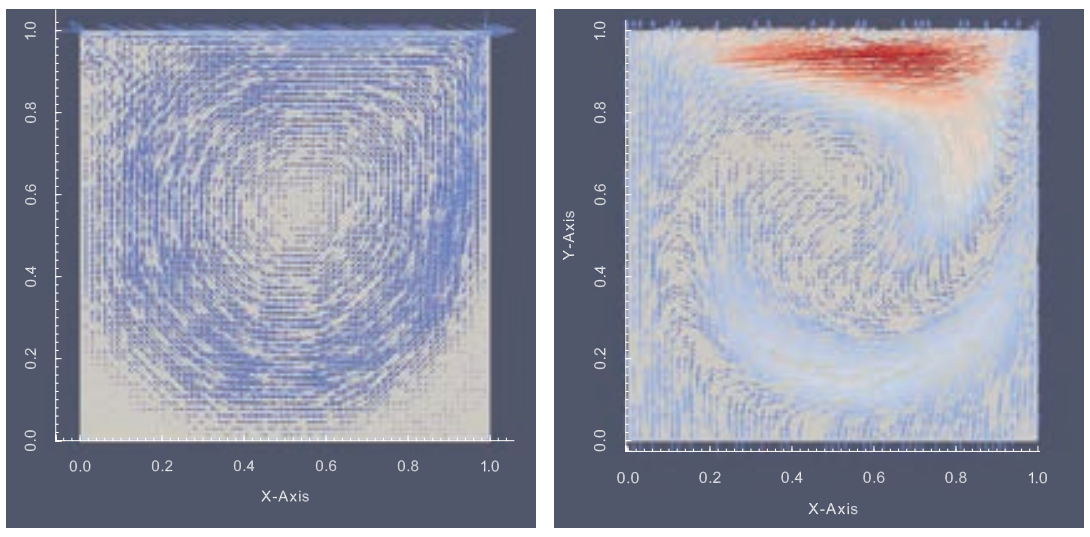

(c) $R_e = 1000, R_m = 200, H_a = 10$.

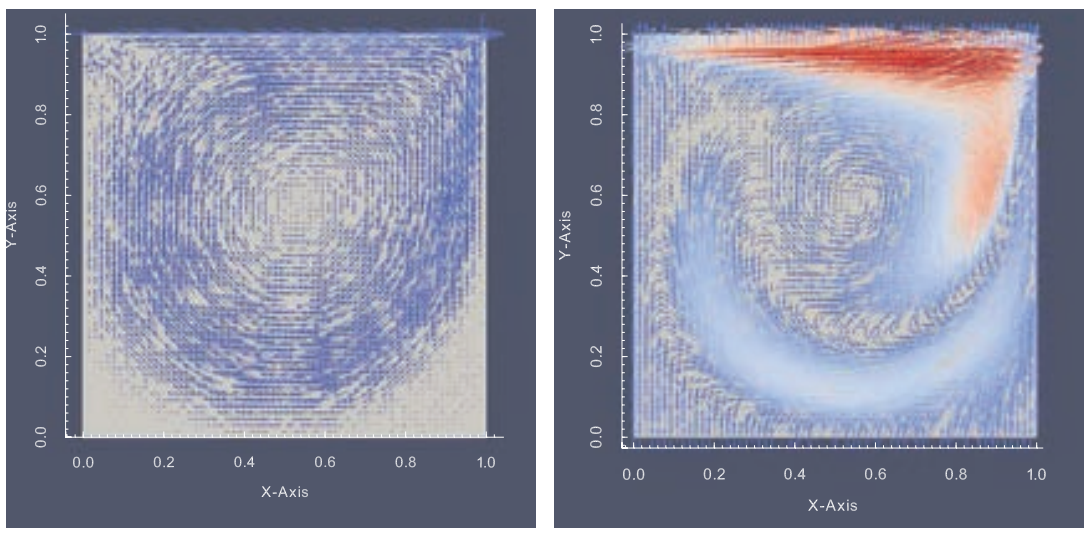

(d) $R_e = 1000, R_m = 400, H_a = 10$.

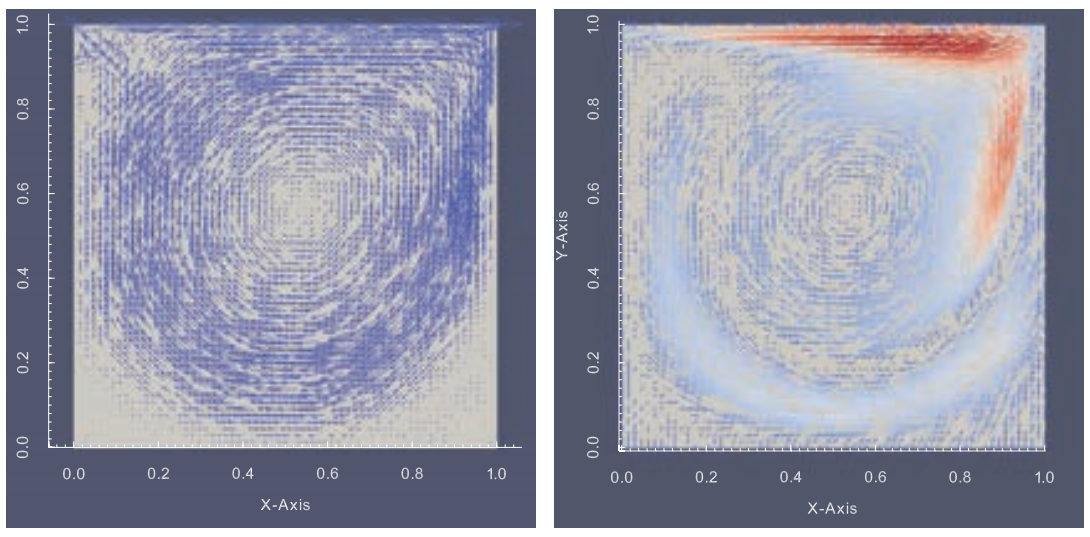

(e) $R_e = 1000, R_m = 1000, H_a = 10$.

FIG. 4. *The steady state of the velocity field* $\boldsymbol{u}$ *and the magnetic field* $\boldsymbol{B}$ *for five set of order parameters, shown in the caption of each subfigure, in which, the left one is the velocity field and the right one is the magnetic field.*

Second, we perform the so-called "magnetic frozen" simulations by varying the magnetic Reynolds number $R_m$. We set the fluid Reynolds number and Hartmann number to be $R_e = 1000, H_a = 10$, and five different values of $R_m = 10, 100, 200, 400, 1000$ in Fig. 4. We observe that the magnetic field bends gradually as $R_m$ grows, which means the convection of magnetic field yields more and more notable influences. When the magnetic field is strong enough, for instance when $R_m = 100, 200, 400, 1000$, the magnetic field form a big eddy at the center of the cavity which is quite similar to the corresponding fluid flow. This feature is called as the magnetic field is *frozen* inside the velocity field. Similar phenomenon had been reported in [29] as well.

**5. Concluding Remarks.** In this paper, we develop an efficient numerical scheme for solving the MHD system. The scheme is (a) fully decoupled, (b) unconditionally energy stable, (c) linear and easy-to-implement. Moreover, we theoretically establish the unconditional energy stability and provide rigorous error estimates for the scheme. A series of numerical simulations, including the convergence test, energy stability test and two physical benchmark problems, are presented to validate the stability and accuracy of the scheme.

**Acknowledgments.** The authors are grateful to Edward G. Phillips and Eric C. Cyr for their kindly help in the numerical simulations.